\documentclass{article}
\usepackage[english
]{babel}
\usepackage[a4paper]{geometry}
\usepackage{amsmath,amssymb,amsthm}

\usepackage{hyperref}
\hypersetup{colorlinks=true,linkcolor=blue,citecolor=teal,
filecolor=magenta,urlcolor=cyan}

\usepackage{tikz}
\usetikzlibrary{calc,decorations.markings}

\usepackage{graphicx}

\newtheorem{theorem}{Theorem}[section]
\newtheorem{lemma}[theorem]{Lemma}
\newtheorem{corollary}[theorem]{Corollary}
\newtheorem{proposition}[theorem]{Proposition}

\newtheorem{conjecture}[theorem]{Conjecture}

\theoremstyle{definition}
\newtheorem{definition}[theorem]{Definition}
\newtheorem*{remark}{Remark}
\newtheorem{example}[theorem]{Example}

\newtheorem*{problem}{Problem}

\newcommand{\BC}{\mathbb{C}}

\newcommand{\BF}{\mathbb{F}}

\newcommand{\BN}{\mathbb{N}}

\newcommand{\cB}{\mathcal{B}}

\newcommand{\cC}{\mathcal{C}}

\newcommand{\cD}{\mathcal{D}}

\newcommand{\cF}{\mathcal{F}}

\newcommand{\cG}{\mathcal{G}}

\newcommand{\cH}{\mathcal{H}}

\newcommand{\cI}{\mathcal{I}}

\newcommand{\cJ}{\mathcal{J}}

\newcommand{\cS}{\mathcal{S}}

\newcommand{\fg}{\mathfrak{g}}

\newcommand{\fh}{\mathfrak{h}}

\newcommand{\fn}{\mathfrak{n}}

\newcommand{\GL}{\mathfrak{gl}}
\newcommand{\Aut}{\mathop{\rm Aut}}

\def\gl{\mathfrak{gl}}
\def\sl{\mathfrak{sl}}

\title{Weight systems and invariants of graphs and embedded graphs}
\author{M.~Kazaryan, S.~Lando\thanks{National Research University Higher School of Economics,
Skolkovo Institute of Science and Technology}}
\date{}

\begin{document}
\maketitle

\tableofcontents
\newpage

\begin{tabular}{cp{0.8\textwidth}}
  $\fg, \langle\cdot,\cdot\rangle$ & 
 Lie algebra endowed with a nondegenerate
 invariant scalar product \\
  $\gl(N)$ & general linear Lie algebra; consists of all $N\times N$-matrices, with commutator as the Lie bracket\\
  $\gl(1|1)$ & general linear $1|1$ Lie superalgebra\\
   $\sl(N)$ & special linear Lie algebra; consists of all
   $N\times N$--matrices having zero trace, with commutator as the Lie bracket\\
  $d$ & dimension of a Lie algebra; in particular $d=N^2$ for $\gl(N)$,  \\
  $C$ & chord diagram \\
  $n$ & number of chords in a chord diagram  \\
  $g(C)$ & intersection graph of a chord diagram $C$\\
  $G$ & abstract simple graph\\
  $\Gamma$ & embedded graph; allowed to have loops and multiple edges\\
  $\cG$ & Hopf algebra of graphs\\
  $\cC$ & Hopf algebra of chord diagrams\\
  $\cB^e$ & Hopf algebra of even delta-matroids\\
  $K_n$ & complete graph on $n$ vertices, as well as the chord diagram with this intersection graph\\
  $K_{n,m}$ & complete bipartite graph with parts of size $n$ and~$m$, as well as the chord diagram with this intersection graph\\
  $\pi$ &  in various Hopf algebras, projection to
  the subspace of primitive elements the kernel of which
  is the subspace of decomposable elements\\
  $C_1,\cdots,C_N$ & Casimir elements in $U(\gl(N))$\\
  $w$ & weight system\\
  $w_\fg$ & weight system associated to a Lie algebra $\fg$\\
  $\bar{w}_\fg$ & $w_\fg(\pi(\cdot))$; composition of a weight system $w_\fg$ and the projection~$\pi$ 
  to the subspace of primitives\\
  $S_m$ & symmetric group\\
  $\sigma$ & permutation\\
  $m$ & number of permuted elements; for example, 
  $m=2n$ for the permutation defined by a chord diagram\\
  $G(\sigma)$ & oriented graph of a permutation~$\sigma$\\
  $\Lambda^*(N)$ & algebra of shifted symmetric polynomials in~$N$ variables\\
  $\phi$ & Harish-Chandra projection\\
  $p_1,\cdots,p_N$ & shifted power sums\\
  $\nu(G)$ & nondegeneracy of a graph~$G$\\
  $A_G(q_1,q_2,\dots)$& Abel polynomial of a graph~$G$\\
   $\chi_G(c)$& chromatic polynomial of a graph~$G$\\
    $S_G(q_1,q_2,\dots)$& Stanley's symmetrized chromatic polynomial of a graph~$G$\\
     $W_G(q_1,q_2,\dots)$& weighted chromatic polynomial of a graph~$G$\\
      $T_C(x,t,s)$& transition polynomial of a chord diagram~$C$\\
       $Q_G(x)$& skew characteristic polynomial of a graph~$G$\\
        $L_G(x)$& interlace polynomial of a graph~$G$\\
\end{tabular}\\

\newpage

The theory of finite type invariants of knots and links
constructed mainly in 1990'ies originating in 
the paper~\cite{V90} leads to a necessity to study weight systems.
A weight system is a function on chord diagrams 
satisfying so-called Vassiliev's $4$-term relations.
Usually one considers weight systems taking values
in a commutative ring. In turn, a chord diagram is
a combinatorial object, a graph of special kind.
One can think of it as of a $3$-regular graph 
together with a distinguished Hamiltonian circuit.
Two such graphs are considered as being isomorphic
if there is an isomorphism taking the Hamiltonian
circuit in the first graph to the one in the second
graph, with orientation preserved.

A $4$-term relation is constructed from a chord diagram
and a pair of its vertices that are neighbors along
the Hamiltonian circuit; it is a linear relation on
the values of a function on four chord diagrams
with one and the same number of vertices.
Since the relations are linear, all weight systems
can be described as solutions to a system of linear equations.
However, both the number of variables (chord diagrams),
and the number of equations ($4$-term relations) grows
very rapidly as the number of vertices grows (thus,
there are more than 600 thousand diagrams with
$2\cdot 9=18$ vertices, and the number of relatons exceeds $5\cdot10^6$). In contrast, the rank of the
equations matrix, that is, the dimension of the space
of chord diagrams modulo $4$-term relations, grows
much slower, and for $9$ chords the rank is only $104$.
There are no universal tools to compute these dimensions,
and only rough upper and lower estimates are known.

Studying spaces of all weight systems must be completed
by constructing specific weight systems and families
of such weight systems. These constructions are
used, in particular, to obtain lower bounds for
the dimensions.
Two main sources of constructions are invariants of
intersection graphs of chord diagrams that satisfy
$4$-term relations for graphs~\cite{L00}, and metrized Lie algebras~\cite{K93,BN95}. 
These two constructions are essentially different.
Graph invariants usually can be easily computed,
but their distinguishing power is rather poor.
In contrast, weight systems associated to Lie algebras
are powerful, but no approaches to their efficient
computation has been known till recently. In the
present paper, we describe both approaches, with 
a stress on recent results.

The problem of computing values of weight systems 
associated with metrized Lie algebras on chord diagrams
is complicated because computations must be done in
a highly noncommutative algebra, which is the universal
enveloping algebra of the Lie algebra.
For the simplest nontrivial case of the Lie algebra
 $\sl(2)$, the computations could be simplified by
 using the Chmutov--Varchenko recurrence relations~\cite{CV97}, 
but they also often prove to be nonsufficient. 
Thus, for long we have not known the values of the
$\sl(2)$-weight system of chord diagrams in which
any two chords intersect one another (that is,
whose intersection graph is a complete graph).
The final answer has been given by P.~Zakorko~\cite{Za22} by proving S.~Lando's 2015 conjecture
about the generating function for these values.
Another recent achievement is the construction of 
recurrence relations for the values of the
$\gl(N)$-weight system~\cite{ZY22}. These relations
are based on an extension of the
$\gl(N)$-weight system to arbitrary permutations
suggested by M.~Kazarian (while chord diagrams are
interpreted as permutations of special kind).

One of the key problems in studying graph invariants
is the one about their extension to embedded graphs.
This question is studied in many recent papers.
In most of them, the invariant of embedded graphs is
defined as the invariant of the underlying abstract 
graph enriched by information about the embedding.
We, however, are interested first of all in extending
weight systems, which are defined by
finite type knot invariants, to weight systems
associated to finite type invariants of links
(a link, in contrast to a knot, has several 
connected components). A chord diagram can be
interpreted as an embedded graph with a single vertex 
on an orientable surface, while links correspond to
generalized chord diagrams, which are embedded graphs, 
whose number of vertices coincides with the number
of link components, on orientable surfaces.

In recent papers~\cite{NZ21,DL22}, an approach 
to extending weight systems and graph invariants
to arbitrary embedded graphs, which is based
on the study of the structure of the corresponding
Hopf algebras. The space of graphs, as well as the
spaces of chord diagrams modulo $4$-term relations
are endowed with natural connected graded 
Hopf algebra structure. We do not know such a structure
on the space spanned by embedded graphs.
It exists, however, on the space of binary 
delta-matroids~\cite{LZ17}, which are combinatorial
objects that accumulate important information about
graphs or chord diagrams, and about embedded graphs.
However, extending weight systems associated to
metrized Lie algebras to embedded graphs with arbitrarily mny vertices remains an open problem.

The paper has the following structure.

In Sec.~\ref{s1}, we recall main definitions
related to chord diagrams, their intersection graphs,
$4$-term relations for chord diagrams and graphs.
We give several key examples of graph invariants
that satisfy $4$-term relations and define thus
weight systems.

In Sec.~\ref{s2}, the Hopf algebra structure on the
space of graphs is described as well as the 
behavior of the graph invariants described above
with respect to this structure. Here we also mention
recent results~\cite{CKL20} relating the Hopf algebra
structure to an integrability property of graph
invariants.

Section~\ref{s3} is devoted to weight systems 
associated to Lie algebras. In addition to known
definitions and results about the $\sl(2)$-weight
system we describe recent results about its values
on chord diagrams whose intersection graphs is
complete bipartite, due to P.~Zinova and M.~Kazarian,
as well as Lando's conjecture (P.~Zakorko's theorem)
about its values on chord diagrams with complete intersection graphs.
This section also describes the extension of the 
$\gl(N)$-weight system to arbitrary permutations
and recurrence relations for its computation.

Section~\ref{s4} is devoted to constructing extensions
of weight systems from chord diagrams to arbitrary
embedded graphs. Here we introduce delta-matroids and
describe the Hopf algebra structure on the space of
even binary delta-matroids.

All necessary definitions and classical results about weight systems can be found in~\cite{CDBook12}.
We reproduce only those definitions and statements 
that are necessary to state recent results.

The authors are partially supported by International Laboratory of Cluster Geometry NRU HSE, RF Government grant, ag. № 075-15-2021-608 dated 08.06.2021.

\section{$4$-term relations}\label{s1}

In 1990, V.~Vassiliev~\cite{V90} introduced the notion of finite type
knot invariant and associated to any knot invariant of order at most~$n$
a function on chord diagrams with~$n$ chords. He showed that any such 
function satisfies $4$-term relations. In 1993, M.~Kontsevich~\cite{K93}
proved the converse: if a function on chord diagrams with~$n$ chords
taking values in a field of characteristic zero satisfies the $4$-term
relations, then it can be obtained by Vassiliev's construction from some
knot invariant of order at most~$n$. (In Vassiliev's theorem, as well as
in Kontsevich's one, the function is subject to an additional
requirement; namely, it must satisfy the so-called one-term relation. 
This requirement, however, vanishes when knots are replaced by framed 
knots, it does not affect the theory, and we will not mention it below).
Functions on chord diagrams satisfying $4$-term relations are called
\emph{weight systems}.

The universal finite type knot invariant constructed by Kontsevich 
(\emph{Kontsevich's integral})
allows one to reconstruct, in principle, a finite type invariant 
from the corresponding weight system. Therefore, studying weight systems
plays a key role in understanding of the nature of finite type knot
invariants. Note, however, that explicit computation of the Kontsevich
integral of a knot is a computationally complicated problem that 
does not have a satisfactory solution up to now.

\subsection{Chord diagrams and $4$-term reltions}

A \emph{chord diagram of order~$n$} is an oriented circle
together with $2n$ pairwise distinct points on it split into $n$
pairs considered up to orientation preserving diffeomorphisms of the circle.
We call the circle carrying the diagram its \emph{Wilson loop}.

Everywhere in the pictures below we assume that the circle is
oriented counterclockwise. The points forming a pair are connected
by a line or arc segment.

A $4$-term relation is defined by a chord diagram and a pair of
chords having neighboring ends in it. We say that 
\emph{a function~$f$ on chord diagrams satisfies Vassiliev's
$4$-term relations} if for any chord diagram~$C$ and any pair of chords
$a,b$ in~$C$ having neighboring ends the equation shown in 
Fig.~\ref{fourtermrelation} holds. 

\begin{figure}[ht]
\[
f\left(\begin{tikzpicture}[baseline={([yshift=-.5ex]current bounding box.center)}]
	\draw[dashed] (0,0) circle (1);
	\filldraw[black] (-.8,0)  node[anchor=west]{$a$};
	\filldraw[black] (.3,0)  node[anchor=west]{$b$};
	\draw[line width=1pt]  ([shift=( 20:1cm)]0,0) arc [start angle= 20, end angle= 70, radius=1];
	\draw[line width=1pt]  ([shift=(110:1cm)]0,0) arc [start angle=110, end angle=160, radius=1];
	\draw[line width=1pt]  ([shift=(250:1cm)]0,0) arc [start angle=250, end angle=290, radius=1];
	\draw[line width=1pt] (45:1) ..  controls (5:0.3) and (-40:0.3)  .. (280:1);
	\draw[line width=1pt] (135:1) ..  controls (175:0.3) and (220:0.3)  .. (260:1);
\end{tikzpicture}\right)  -
f\left(\begin{tikzpicture}[baseline={([yshift=-.5ex]current bounding box.center)}]
	\filldraw[black] (-.6,0)  node[anchor=west]{$a$};
	\filldraw[black] (.2,0)  node[anchor=west]{$b$};
	\draw[dashed] (0,0) circle (1);
	\draw[line width=1pt]  ([shift=( 20:1cm)]0,0) arc [start angle= 20, end angle= 70, radius=1];
	\draw[line width=1pt]  ([shift=(110:1cm)]0,0) arc [start angle=110, end angle=160, radius=1];
	\draw[line width=1pt]  ([shift=(250:1cm)]0,0) arc [start angle=250, end angle=290, radius=1];
	\draw[line width=1pt] (45:1) ..  controls (-5:0.1) and (-50:0.1)  .. (260:1);
	\draw[line width=1pt] (135:1) ..  controls (185:0.1) and (225:0.1)  .. (280:1);
\end{tikzpicture}\right)  =
f\left(\begin{tikzpicture}[baseline={([yshift=-.5ex]current bounding box.center)}]
	\filldraw[black] (-.4,.4)  node[anchor=west]{$a$};
	\filldraw[black] (.4,0)  node[anchor=west]{$b$};
	\draw[dashed] (0,0) circle (1);
	\draw[line width=1pt]  ([shift=( 20:1cm)]0,0) arc [start angle= 20, end angle= 70, radius=1];
	\draw[line width=1pt]  ([shift=(110:1cm)]0,0) arc [start angle=110, end angle=160, radius=1];
	\draw[line width=1pt]  ([shift=(250:1cm)]0,0) arc [start angle=250, end angle=290, radius=1];
	\draw[line width=1pt] (35:1) ..  controls (0:0.3) and (-45:0.3)  .. (280:1);
	\draw[line width=1pt] (135:1) ..  controls (105:0.5) and (85:0.5)  .. (55:1);
\end{tikzpicture}\right)  -
f\left(\begin{tikzpicture}[baseline={([yshift=-.5ex]current bounding box.center)}]
    \filldraw[black] (-.4,.3)  node[anchor=west]{$a$};
	\filldraw[black] (.2,0)  node[anchor=west]{$b$};
	\draw[dashed] (0,0) circle (1);
	\draw[line width=1pt]  ([shift=( 20:1cm)]0,0) arc [start angle= 20, end angle= 70, radius=1];
	\draw[line width=1pt]  ([shift=(110:1cm)]0,0) arc [start angle=110, end angle=160, radius=1];
	\draw[line width=1pt]  ([shift=(250:1cm)]0,0) arc [start angle=250, end angle=290, radius=1];
	\draw[line width=1pt] (55:1) ..  controls (5:0.1) and (-40:0.1)  .. (270:1);
	\draw[line width=1pt] (135:1) ..  controls (105:0.4) and (65:0.4)  .. (35:1);
\end{tikzpicture}\right)
\]
\caption{$4$-term relation for chord diagrams}
\label{fourtermrelation}
\end{figure}

The equation in the picture has the following meaning. All the four
chord diagrams entering it can have an additional tuple of chords
whose ends belong to the dashed arcs, which is one and the same
in all the four diagrams. Out of the four ends of the two distinguished
chords $a,b$, both end of~$b$ as well as one end of~$a$ are fixed,
while the second end of~$a$ takes successively all the four positions
close to the ends of~$b$. The specific position of the second end of~$a$
determines the chord diagram in the relation.
Note that it does not matter whether the two chords~$a$ and~$b$
intersect one another or not: the case of intersecting chord reduces to
that of nonintersecting ones by multiplying the relation shown in
Fig.~\ref{fourtermrelation} by $-1$.

Below, we write down the $4$-term relation in the form
\begin{eqnarray}
f(C)-f(C'_{ab})=f(\widetilde{C}_{ab})-f(\widetilde{C}'_{ab})
\end{eqnarray}
and say that the chord diagram~$C'_{ab}$ is obtained from~$C$
by \emph{Vassiliev's first move} on the pair of chords $a,b$,
the chord diagram 
$\widetilde{C}_{ab}$ is the result of \emph{Vassiliev's second move},
and $\widetilde{C}'_{ab}$ is the result of the composition of the
two moves; note that the first and the second move, when elaborated
on one and the same pair of chords, commute with one another,
and the order in which they are performed is irrelevant.

The definition of the $4$-term relation uses subtraction, and
we are planning to multiply chord diagrams in what follows.
That is why we will usually consider weight systems with values
in a commutative ring, making the ring explicit in specific examples.

\subsection{Intersection graphs}

\begin{figure}[h]
\begin{center}
\begin{picture}(200,100)(60,0)
\thicklines
\multiput(12,72)(90,0){4}{\circle{40}}
\multiput(55,72)(180,0){2}{{\scriptsize $-$}}
\put(145,72){$=$}
\multiput(5,91)(90,0){4}{\line(-1,-3){10}}
\multiput(0,94)(90,0){4}{{\scriptsize E}}
\multiput(-3,84)(90,0){4}{\line(1,0){31}}
\multiput(-9,87)(90,0){4}{{\scriptsize D}}
\multiput(-8,72)(90,0){4}{\line(3,-1){36}}
\multiput(-14,73)(90,0){4}{{\scriptsize F}}
\multiput(32,72)(90,0){4}{\line(-5,-2){35}}
\multiput(34,73)(90,0){4}{{\scriptsize B}}
\multiput(12,92)(90,0){4}{\line(1,-5){7.5}}
\multiput(10,95)(90,0){4}{{\scriptsize C}}
\put(18,91){\line(-1,-5){7.6}}
\multiput(18,94)(90,0){4}{{\scriptsize A}}
\put(108,91){\line(-2,-3){20.9}}
\put(198,91){\line(1,-2){13}}
\put(288,91){\line(1,-1){12.3}}

\multiput(2,30)(90,0){4}{\circle*{3}}
\multiput(22,30)(90,0){4}{\circle*{3}}
\multiput(2,-10)(90,0){4}{\circle*{3}}
\multiput(22,-10)(90,0){4}{\circle*{3}}
\multiput(-8,10)(90,0){4}{\circle*{3}}
\multiput(32,10)(90,0){4}{\circle*{3}}
\multiput(55,10)(180,0){2}{{\scriptsize $-$}}
\put(145,10){$=$}
\multiput(-6,-16)(90,0){4}{{\scriptsize A}}
\multiput(24,-16)(90,0){4}{{\scriptsize B}}
\multiput(-6,32)(90,0){4}{{\scriptsize E}}
\multiput(22,32)(90,0){4}{{\scriptsize D}}
\multiput(-14,12)(90,0){4}{{\scriptsize F}}
\multiput(34,12)(90,0){4}{{\scriptsize C}}

\multiput(-8,10)(90,0){4}{\line(1,0){39}}
\multiput(-8,10)(90,0){4}{\line(1,2){10}}
\multiput(-8,10)(90,0){2}{\line(1,-2){10}}
\multiput(-8,10)(90,0){4}{\line(3,-2){30}}
\multiput(32,10)(90,0){2}{\line(-3,-2){30}}
\multiput(32,10)(90,0){4}{\line(-1,2){10}}
\multiput(32,10)(90,0){4}{\line(-1,-2){10}}
\multiput(2,30)(90,0){4}{\line(1,0){20}}
\multiput(2,-10)(180,0){2}{\line(1,0){20}}
\multiput(2,-10)(90,0){4}{\line(1,2){20}}

\end{picture}
\end{center}
		\caption{A $4$-term relation for chord diagrams and
		corresponding intersection graphs}
		\label{f4tig}
\end{figure}
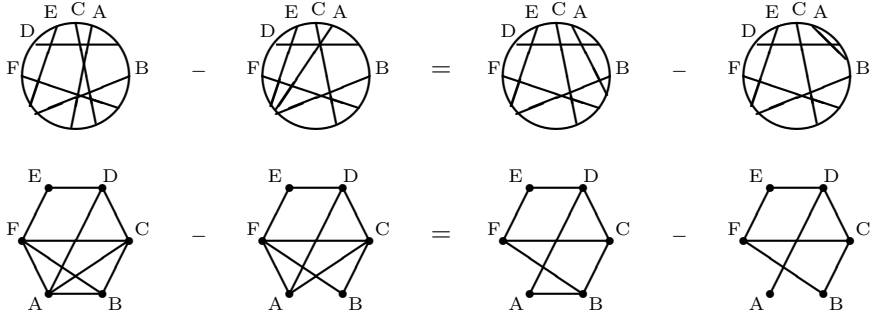

Invariants of intersection graphs of chord diagrams serve
as one of the main sources of weight systems,

The \emph{intersection graph $g(C)$ of a chord diagram~$C$} 
is the simple graph the set of whose vertices is in one-to-one 
correspondence with the set of chords of~$C$, and two vertices
are connected by an edge iff the corresponding chords intersect.
(Two chords in a chord diagram intersect one another if their
ends follow the circle in alternating order; the intersection
point of chords in pictures is not a vertex of the chord diagram.) 
Fig.~\ref{f4tig} shows an example of a $4$-term relation for 
a chord diagram and its intersection graph. Only the chord
diagrams and their intersection graphs are depicted; there is no
reference to the value of the function.

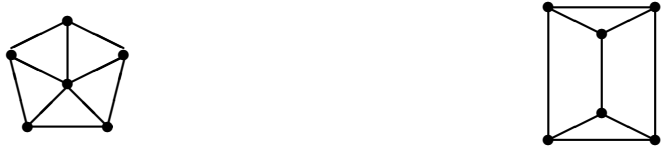
\begin{figure}[h]

\begin{center}
\begin{picture}(300,60)(10,40)

\thicklines

\put(55,51){\circle*{4}}
\put(55,75){\circle*{4}}
\put(34,62){\circle*{4}}
\put(76,62){\circle*{4}}
\put(40,35){\circle*{4}}
\put(70,35){\circle*{4}}
\put(55,50){\line(0,1){24}}
\put(55,51){\line(2,1){21}}
\put(55,51){\line(-2,1){21}}
\put(55,50){\line(1,-1){16}}
\put(55,50){\line(-1,-1){16}}
\put(40,35){\line(1,0){32}}
\put(55,75){\line(2,-1){21}}
\put(55,75){\line(-2,-1){21}}
\put(40,35){\line(-1,4){7}}
\put(70,35){\line(1,4){7}}

\put(255,70){\circle*{4}}
\put(255,40){\circle*{4}}
\put(235,80){\circle*{4}}
\put(275,80){\circle*{4}}
\put(235,30){\circle*{4}}
\put(275,30){\circle*{4}}
\put(235,30){\line(0,1){50}}
\put(235,30){\line(1,0){40}}
\put(235,80){\line(1,0){40}}
\put(275,30){\line(0,1){50}}
\put(255,40){\line(0,1){30}}
\put(235,30){\line(2,1){20}}
\put(235,80){\line(2,-1){20}}
\put(275,80){\line(-2,-1){20}}
\put(275,30){\line(-2,1){20}}

\end{picture}
\end{center}

		\caption{$5$-wheel and $3$-prism, the two graphs with~$6$
		vertices that are not intersection graphs}
		\label{fR3}
\end{figure}

Any graph with at most~$5$ vertices is the intersection graph
of some chord diagram. There are two graphs with $6$ vertices that
are not intersection graphs; they are shown in Fig.~\ref{fR3}. 
As~$n$ grows, the fraction of intersection graphs among simple graphs with~$n$ vertices grows rapidly.
A variety of complete sets of obstacles for a graph to be an intersection graph is known.

On the other hand, certain intersection graphs can be realized
as intersection graphs of several, not pairwise isomorphic 
chord diagrams. Consider, for example, the chain on $n$ vertices,

\begin{center}
  \begin{picture}(80,40)(-50,-20)
     \unitlength=20pt
  \put(0,0){\makebox(0,0){$
              \underbrace{
              \makebox(5.4,0){
                \begin{picture}(5,0.2)(0,-0.2)
                \put(0,0){\circle*{0.15}}
                \put(0,0){\line(1,0){1}}
                \put(1,0){\circle*{0.15}}
                \put(1,0){\line(1,0){1}}
                \put(2,0){\circle*{0.15}}
                \put(2,0){\line(1,0){0.5}}
                \put(3.08,0){\makebox(0,0){\ldots}}
                \put(3.5,0){\line(1,0){0.5}}
                \put(4,0){\circle*{0.15}}
                \put(4,0){\line(1,0){1}}
                \put(5,0){\circle*{0.15}}
                \end{picture}
                           }
                         }_{\scriptstyle n\ {\rm vertices}}$
                         }
           }
  \end{picture}
\end{center}
For $n=5$, there are three distinct chord diagrams with this intersection graph:
%
%
\def\celvn
{  \bezier{40}(-1,0)(-1,0.577)(-0.5,0.866)
   \bezier{40}(-0.5,0.866)(0,1.155)(0.5,0.866)
   \bezier{40}(0.5,0.866)(1,0.577)(1,0)
   \bezier{40}(1,0)(1,-0.577)(0.5,-0.866)
   \bezier{40}(0.5,-0.866)(0,-1.155)(-0.5,-0.866)
   \bezier{40}(-0.5,-0.866)(-1,-0.577)(-1,0)
}
\def\vert
{   \put(0,1){\circle*{0.15}}
    \put(0,-1){\circle*{0.15}}
    \put(0,1){\line(0,-1){2}}
}
%
%
\begin{center}
\begin{picture}(180,60)(0,-30)
  \setlength{\unitlength}{20pt}
  \put(0,0){
     \begin{picture}(0,0)
     \celvn\vert
     \put(0.866,0.5){\circle*{0.15}}
     \put(-0.5,0.866){\circle*{0.15}}
     \bezier{40}(-0.5,0.866)(0.086,0.322)(0.866,0.5)
     \put(0.5,-0.866){\circle*{0.15}}
     \put(-0.866,-0.5){\circle*{0.15}}
     \bezier{40}(-0.866,-0.5)(-0.086,-0.322)(0.5,-0.866)
     \put(0.5,0.866){\circle*{0.15}}
     \put(1,0){\circle*{0.15}}
     \bezier{40}(1,0)(0.433,0.25)(0.5,0.866)
     \put(-0.5,-0.866){\circle*{0.15}}
     \put(-1,0){\circle*{0.15}}
     \bezier{40}(-1,0)(-0.433,-0.25)(-0.5,-0.866)
     \end{picture}
           }
  \put(4,0){
     \begin{picture}(0,0)
     \celvn
     \put(-0.5,0.866){\circle*{0.15}}
     \put(-0.5,-0.866){\circle*{0.15}}
     \put(-0.5,-0.866){\line(0,1){1.732}}
     \put(-0.866,0.5){\circle*{0.15}}
     \put(0.5,0.866){\circle*{0.15}}
     \bezier{40}(-0.866,0.5)(-0.086,0.322)(0.5,0.866)
     \put(0,1){\circle*{0.15}}
     \put(0.866,0.5){\circle*{0.15}}
     \bezier{40}(0.866,0.5)(0.25,0.433)(0,1)
     \put(-0.866,-0.5){\circle*{0.15}}
     \put(0.5,-0.866){\circle*{0.15}}
     \bezier{40}(-0.866,-0.5)(-0.086,-0.322)(0.5,-0.866)
     \put(0,-1){\circle*{0.15}}
     \put(0.866,-0.5){\circle*{0.15}}
     \bezier{40}(0.866,-0.5)(0.25,-0.433)(0,-1)
     \end{picture}
           }
  \put(8,0){
     \begin{picture}(0,0)
     \celvn\vert
     \put(-0.866,0.5){\circle*{0.15}}
     \put(0.5,0.866){\circle*{0.15}}
     \bezier{40}(0.5,0.866)(-0.086,0.322)(-0.866,0.5)
     \put(-0.5,-0.866){\circle*{0.15}}
     \put(0.866,-0.5){\circle*{0.15}}
     \bezier{40}(0.866,-0.5)(0.086,-0.322)(-0.5,-0.866)
     \put(-0.5,0.866){\circle*{0.15}}
     \put(-1,0){\circle*{0.15}}
     \bezier{40}(-1,0)(-0.433,0.25)(-0.5,0.866)
     \put(0.5,-0.866){\circle*{0.15}}
     \put(1,0){\circle*{0.15}}
     \bezier{40}(1,0)(0.433,-0.25)(0.5,-0.866)
     \end{picture}
           }
\end{picture}
\end{center}

\subsection{$4$-term relations for graphs}\label{ss4tg}

$4$-term relations for graphs were introduced in~\cite{L00}.
We say that a graph invariant~$f$ 
\emph{satisfies the $4$-term relation}
if for any graph~$G$ and any pair $a,b$ of its vertices we have
\begin{eqnarray}
f(G)-f(G'_{ab})=f(\widetilde{G}_{ab})-f(\widetilde{G}'_{ab});
\end{eqnarray}
here the graph~$G'_{ab}$ is obtained from~$G$ 
by switching adjacency between the vertices~$a$ and~$b$,
the graph $\widetilde{G}_{ab}$ is obtained from~$G$ by switching
adjacency of~$a$ to all the vertices in~$G$ adjacent to~$b$,
and the graph $\widetilde{G}'_{ab}$ is the result of composition of
these two operations.

\begin{remark}
The transformation $G\mapsto G'_{ab}$ is symmetric with
respect to~$a$ and~$b$,
$G'_{ab}=G'_{ba}$. In turn, the second Vassiliev move for graphs
 $G\mapsto \widetilde{G}_{ab}$ is not symmetric, and
the graph $\widetilde{G}_{ba}$is not, as a rule,
isomorphic to $\widetilde{G}_{ba}$.
\end{remark}

Graph invariants satisfying $4$-term relations are called
\emph{$4$-invariants}. Similarly to the case of weight systems,
we assume that $4$-invariants take values in a commutative ring.

A direct study of the way the intersection graph of a chord diagram
changes under Vassiliev moves proves the following statement.

\begin{theorem}[\cite{L00}]
If~$f$ is a $4$-invarint of graphs, then
$f\circ g$ is a weight system on chord diagrams.
\end{theorem}

An example of a $4$-term relation is given in Fig.~\ref{f4tig}.
Note, however, that even if a graph~$G$ is an intersection graph,
this does not mean that the same is true for the other three
graph in a $4$-term relation.

Therefore, any $4$-invariant determines a weight system,
and hence, by Kontsevich's theorem, a knot invariant.
Below, we give several examples of $4$-invariants.

The following question posed by S.~Lando remains open for
the last two decades:

\begin{problem}
Is it true that modulo $4$-term relations for graphs any graph
is a linear combination of intersection graphs?
\end{problem}

Computer computations (I.~A.~Dynnikov, private communication)
show that this is true for all graphs with at most~$8$ vertices.

\subsubsection{Chromatic polynomial}
One of historically first examples of $4$-invariants 
of simple graphs is the chromatic polynomials.
Let~$G$ be a graph, and let $V(G)$ be the set of its vertices.
The 
\emph{chromatic polynomial} $\chi_G(c)$
of~$G$ is the number of regular colorings of the set
$V(G)$ in~$c$ colors, that is, the number of mappings
$V(G)\to\{1,2,\dots,c\}$ such that for any two adjacent
(connected by an edge) vertices the values of the mapping are different.

It is well known that the chromatic polynomial satisfies
the \emph{contraction--deletion relation}:
\begin{equation}
\chi_G(c)=\chi_{G'_e}(c)-\chi_{G''_e}(c)
\end{equation}
for any graph~$G$ and any edge~$e$ in it;
here the graph $G'_e$ is obtained from~$G$ be deleting the edge~$e$,
and  $G''_e$ is the result of contracting this edge.
As an edge is contracted, it is deleted, its two ends become a single
new vertex of the graph, and any multiple edge that can appear after
that is replaced by a single edge.

The contraction--deletion relation can be proved easily: 
it corresponds to splitting of the set of regular colorings of the 
vertices of the graph $G'_e$ in~$c$ colors into two disjoint subsets,
namely, those colorings where the ends of~$e$ are colored differently
(these are exactly the regular colorings of the vertices of~$G$)
and those where these ends are colored with the same color
(the latter correspond one-to-one to regular colorings
of the vertices of $G''_e$). Since the chromatic polynomial 
of a discrete (edge-free) graph on~$n$ vertices is~$c^n$, 
the contraction-deletion relation proves that, for a graph~$G$
on~$n$ vertices, $\chi_G(c)$ indeed is a polynomial of degree~$n$
having the leading coefficient~$1$.

\begin{theorem}[\cite{CDL94}]
Chromatic polynomial is a $4$-invariant,
\end{theorem}

Indeed, applying the contraction--deletion relation
to a graph~$G$ and an edge~$e$ having the ends $a,b$
in it, we conclude that
$$
\chi_{G}(c)-\chi_{G'_{ab}}(c)=-\chi_{G''_{ab}}(c).
$$
In turn,
$$
\chi_{\widetilde{G}_{ab}}(c)-\chi_{\widetilde{G}'_{ab}}(c)
=-\chi_{\widetilde{G}''_{ab}}(c).
$$
Verifying that the natural identification of the sets of vertices
of the two graphs $G''_{ab}$ and $\widetilde{G}''_{ab}$
establishes their isomorphism completes the proof.

\subsubsection{Weighted chromatic polynomial
 \\ (Stanley's symmetrized chromatic polynomial)}

Chromatic polynomial is an important graph invariant;
however, it is not distinguishingly powerful.
For example, the chromatic polynomial of any tree on~$n$
vertices is $c(c-1)^{n-1}$, while the number of pairwise nonisomorphic
trees on~$n$ vertices grows very rapidly as~$n$ grows.
Weighted chromatic polynomial, more widely known under the name
of Stanley's symmetrized chromatic polynomial, is a much finer
graph invariant. 

In order to define weighted chromatic polynomial, we will need
the notion of weighted graph.

\begin{definition}
A \emph{weighted graph} is a simple graph~$G$ together with
a positive integer, a \emph{weight},  assigned to each of its vertices.
The \emph{weight} of a weighted graph is the sum of the weights of
its vertices.
\end{definition}

\begin{definition}[\cite{CDL94}]
The \emph{weighted chromatic polynomial} 
is the weighted graph invariant $G\mapsto W_G(q_1,q_2,\dots)$ 
taking values in the ring of polynomials in infinitely many variables
$q_1,q_2,\dots$ and possessing the following properties:
\begin{itemize}
    \item the weighted chromatic polynomial of the graph on one
    vertex of weight~$n$ is~$q_n$;
    \item the weighted chromatic polynomial is multiplicative, $W_{G_1\sqcup G_2}=W_{G_1}W_{G_2}$
    for a disjoint union of arbitrary weighted graphs $G_1,G_2$;
    \item the weighted chromatic polynomial satisfies the weighted
    contraction--deletion relation
    \begin{equation}
        W_G=W_{G'_e}+W_{G''_e},
    \end{equation}
    where $e$ is an edge of the weighted graph~$G$,
    $G'_e$ denotes the result of deleting an edge~$e$ from~$G$,
    and $G''_e$ is the result of contracting~$e$;
    contraction is defined in the same way as for simple graphs,
    and the new vertex is assigned the weight equal to the sum of 
    the weights of the two ends of~$e$.
\end{itemize}
\end{definition}

\begin{theorem}[\cite{CDL94}]\label{tHawg}
There is a unique weighted graph invariant possessing the
above properties. The weighted chromatic polynomial $W_G$
of a weighted graph~$G$ of weight~$n$ is a quasihomogeneous 
polynomial of the variables~$q_1,q_2,\dots$ of degree~$n$,
if we set the weight of variables~$q_k$ equal to~$k$, for $k=1,2,\dots$.
\end{theorem}

A simple graph can be considered as a weighted graph, with the
weight of each of its vertices set to be~$1$. Thus, any weighted
graphs invariant defines therefore a simple graphs invariant.

\begin{theorem}
The weighted chromatic polynomial is a $4$-invariant.
\end{theorem}

Indeed, similarly to the proof for chromatic polynomial,
we remark that the natural identification of the set
of vertices of the graphs $G''_{ab}$ and $\widetilde{G}''_{ab}$ 
establishes an isomorphism of these graphs as weighted graphs.

In 1995, R.~Stanley introduced the notion of symmetrized chromatic
polynomial.

\begin{definition}[\cite{S95}]
Let~$G$ be a simple graph. A \emph{coloring} of the set of vertices $V(G)$ of~$G$
into an infinite set of colors is a mapping
$\beta:V(G)\to X=\{x_1,x_2,\dots\}$. To each coloring~$\beta$,
one associates a monomial $x_\beta$ of degree $n=|V(G)|$ in the variables $x_1,x_2,\dots$, which is equal to the product of 
the values of~$\beta$ on all the vertices of~$G$.
The \emph{symmetrized chromatic polynomial}
$S_G(x_1,x_2,\dots)$ of~$G$ is the infinite sum
$$
S_G(x_1,x_2,\dots)=
\sum_{\substack{\beta:V(G)\to X\\ \beta\text{ regular}}}
x_\beta,
$$
where the coloring~$\beta$ is said to be \emph{regular}
if it takes any two adjacent vertices to different elements of~$X$.
\end{definition}

By obvious reasons, the symmetrized chromatic polynomial is a
symmetric function in the variables~$X$. Stanley's conjecture 
asserts that the symmetrized chromatic polynomial distinguishes 
between any two trees. It is confirmed for trees with up to~$29$ vertices~\cite{HJ19}, indicating thus that the symmetrized chromatic
polynomial is a much finer graph invariant than the ordinary
chromatic polynomial $\chi_G$.

Being a symmetric polynomial in the variables~$X$,
any symmetrized chromatic polynomial can be expressed as a
polynomial in any basis in the space of symmetric polynomials.
In particular, one can choose the power sums
$$
p_k=\sum_{i=1}^\infty x_i^k
$$
for this basis. When expressed in this form, the symmetrized
chromatic polynomial $S_G$ becomes a finite quasihomogeneous 
polynomial of degree~$n=|V(G)|$ in the variables $p_1,p_2,\dots$
if we set the degree of the variable~$p_k$ equal to~$k$, for
$k=1,2,\dots$. The following statement establishes a relationship
between Stanley's symemtrized chromatic polynomial and
the weighted chromatic polynomial.

\begin{theorem}[\cite{NW99}]
Stanley's symmetrized chromatic polynomial, when expressed
in terms of power sums, under the substitution
 $p_k=(-1)^{n-k}q_k$, $k=1,2,\dots$, becomes the weighted chromatic
 polynomial; here $n=|V(G)|$.
\end{theorem}

As a corollary, Stanley's symmetrized chromatic polynomial
is a $4$-invariant of graphs.

In turn, the ordinary chromatic polynomial is a specialization
of Stanley's symmetrized one: the latter transforms into
the former under the substitution $p_i=c$, for $i=1,2,3,\dots$.

\subsubsection{Interlace polynomial} 

Initially, the interlace polynomial has been defined
as a function on oriented graphs having two-in and two-out edges at each vertex. The definition appeared
first in the paper~\cite{ABS04}, which is devoted to
DNA sequencing. Later, the definition has been extended
to arbitrary simple graphs. We start with defining the
pivot operation.

Let $G$ be a simple graph. For any pair
$a,b$ of adjacent vertices of~$G$, all the other
its vertices are split into four classes, namely:

\begin{enumerate}
    \item the vertices adjacent to~$a$, and not to~$b$;
    \item the vertices adjacent to~$b$, and not to~$a$;
    \item the vertices adjacent to both~$a$, and~$b$;
    \item the vertices adjacent neither to~$a$, nor to~$b$.
\end{enumerate}

The \emph{pivot} $G^{ab}$ of~$G$ along the edge~$ab$ 
is the graph obtained by switching the adjacency between
any two vertices in the first three classes iff they
belong to different classes.

The \emph{interlace polynomial} $L_G(x)$ is defined 
by the following recurrence relations:
\begin{enumerate}
    \item if~$G$ does not have edges, then $L_G(x)=x^n$,
    where $n$  is the number of vertices in~$G$;
    \item for any edge $ab$ in~$G$, we have
    $$
    L_G(x)=L_{G\setminus a}(x)+L_{G^{ab}\setminus b}(x),
    $$
\end{enumerate}
where $G\setminus a$ denotes the graph obtained from~$G$
by deleting the vertex~$a$ and all the edges incident
to it.

In~\cite{ABS04}, it is proved that the interlace polynomial
is well defined: the result of its calculation is 
independent of the order in which the pivots are
applied. The definition immediately implies that 
$L_G(x)$ indeed is a polynomial in~$x$
whose degree is $n=|V(G)|$.

\begin{theorem}
The interlace polynomial satisfies $4$-term relations.
\end{theorem}

This theorem is proved in~\cite{NN} and,  in a different way, in~\cite{K20}.

\subsubsection{Transition polynomial}
Transition polynomial for $4$-regular graphs endowed
with an oriented Eulerian circuit was intro\-duced in~\cite{J90}.
By contracting each chord of a chord diagram~$C$ we
make the latter into a $4$-regular graph in which 
the Wilson loop of the chord diagram becomes an
oriented Eulerian circuit. The transition polynomial
of this graph defines, therefore, a mapping from the
set of chord diagrams to a space of polynomials.

Conversely, to a $4$-regular graph with a distinguished oriented Eulerian circuit a chord diagram is associated;
the number of chords in it coincides with the number of
vertices in the graph. The Eulerian circuit turns into 
the Wilson loop of the chord diagram, and the vertices 
of the graph become the chords.

Here we define the transition polynomial for a chord diagram.
In order to define it, we will need the notion of transition.
Let~$C$ be a chord diagram. Each chord in~$C$ can be replaced by
a ribbon in one of the following three ways; we encode these ways
by the Greek letters $\chi$, $\phi$, or $\psi$:
\begin{itemize}
    \item $\phi$ if this an ordinary ribbon;
    \item $\chi$ if this is a half-twisted ribbon;
    \item $\psi$ if this is no ribbon at all, that is, if we
    simply erase the chord.
\end{itemize}
A choice of a transition for each ribbon, i.e.,
a mapping  $V(C)\to\{\phi,\chi,\psi\}$ taking the set of chords
$V(C)$ to the set of transition types, is called
a \emph{state} of the chord diagram~$D$.

Choose a \emph{weight function} 
$w:\{\phi,\chi,\psi\}\to K$ that associates to each of the three
Greek letters an element of a commutative ring~$K$.

The \emph{weighted transition polynomial} of a chord diagram~$C$
is the polynomial
$$
T_C(x)=\sum_{s:V(C)\to\{\phi,\chi,\psi\}}
\prod_{v\in V(C)} w(s(v))x^{c(s)-1},
$$
where the summation is carried over all states of~$C$,
and  $c(s)$ denotes the number of connected components of
the boundary of the surface obtained from~$C$ by attaching
a disk to the Wilson loop with the chords replaced by
the corresponding ribbons.

\begin{theorem}[\cite{DZ22}]
By choosing $w(\chi)=t$,
$w(\phi)=-t$, $w(\psi)=s$, for the weight function, we make 
the transition polynomial into a weight system taking values
in the ring of polynomials $\BC[t,s,x]$.
\end{theorem}

\subsubsection{Skew characteristic polynomial}\label{sssSCP} 

Let $G$ be a simple graph with $n=|V(G)|$ vertices.
Number the vertices by the numbers from~$1$ to~$n$
in an arbitrary way and associate to the chosen numbering the
adjacency matrix of~$G$.
The \emph{adjacency matrix} $A(G)$ is an
$n\times n$-matrix over the two-element field~$\BF_2$,
containing~$1$ in the cell $(i,j)$
if the vertices numbered~$i$ and~$j$ are connected by an edge,
and containing~$0$ otherwise. In particular, the matrix $A(G)$
is symmetric and has zeroes on the diagonal. The isomorphism
class of~$G$ is reconstructed uniquely from its adjacency matrix.
Various characteristics of the adjacency matrix, for example, its characteristic polynomial, play a key role in the study of graphs.

\begin{definition}
A graph~$G$ is \emph{nondegenerate} (respectively, 
\emph{degenerate}) if its adjacency matrix $A(G)$ is nondegenerate
(respectively, \emph{degenerate}) over~$\BF_2$.
The \emph{nondegeneracy} of a graph is the graph invariant taking
values in the field~$\BC$ and equal to~$1$ for a nondegenerate graph
and to~$0$ for a degenerate one.
\end{definition}

The adjacency matrix of any graph with an odd number of veertices
is degenerate since it is antisymmetric over the field $\BF_2$,
whence the nondegeneracy of any graph with an odd number of vertices
is~$0$. A graph with an even number of vertices can be either
nondegenerate or degenerate. We denote the nondegeneracy of a
graph~$G$ by $\nu(G)$.

\begin{lemma}\label{lnd}
Nondegeneracy is invariant with respect to the second Vassiliev move.
\end{lemma}

Indeed, let $G$ be  graph, and let $a,b$ be its vertices.
Pick a numbering of the vertices of~$G$ such that~$a$ is assigned number~$1$
and~$b$ is assigned number~$2$. Then the adjacency matrix  $A(\widetilde{G}_{ab})$ of
$\widetilde{G}_{ab}$ is
\begin{eqnarray}\label{eVsm}
A(\widetilde{G}_{ab})=I_{12}^tA(G)I_{12},
\end{eqnarray}
where $I_{12}$ is the square $|V(G)|\times |V(G)|$-matrix over~$\BF_2$,
which coincides with the identity matrix everywhere with the exception
of the cell $(1,2)$ whose value equals~$1$,
and $I_{12}^t$ is its transpose matrix. The assertion of the lemma
follows from the fact that both matrices $I_{12}$ and $I_{12}^t$
are nondegenerate.

\begin{definition}[\cite{DL22}]
The \emph{skew characteristic polynomial} of a graph~$G$ is the polynomial 
$$
Q_G(u)=\sum_{U\subset V(G)} \nu(G|_U)u^{|V(G)|-|U|},
$$
where the summation is carried over all subsets~$U$ of the set $V(G)$
of vertices of the graph, and $G|_U$ denotes the subgraph of~$G$
induced by a subset~$U$ of its vertices.
\end{definition}

The skew characteristic polynomial of a graph with an even number of vertices
is an even function of its argument, and if the number of vertices is
odd, then the skew characteristic polynomial is odd as well.

Lemma~\ref{lnd} implies 

\begin{proposition}
The skew characteristic polynomial is a $4$-invariant of graphs.
\end{proposition}

Let us explain why the skew characteristic polynomial carries this name.
Let $C$ be a chord diagram, and let  $g(C)$ be its intersection graph.
Pick the following orientation of the graph. Choosing a point
$\alpha$ on the Wilson loop that is not an end of any chord,
cut the circle at this point and develop it into a horizontal line
whose orientation inherits that of the circle. Under this transformation,
the chords become half circles (arcs), and two vertices in $g(C)$
are connected by an edge iff the corresponding arcs intersect
one another. Orient each edge in  $g(C)$ from the vertex corresponding
to the arc whose left end precedes the left end of the arc corresponding
to the second vertex of the edge. If we number the vertices of the graph
in the order of the left ends of the corresponding arcs,
then each edge is oriented from the vertex with a smaller number to that
with a greater one. The orientation of the intersection graph obtained
in this way depends on the cut point~$\alpha$ of the circle;
denote the corresponding directed graph by
$\vec{g}_\alpha(C)$.

The intersection matrix of a directed graph with~$n$ numbered 
vertices is the integer-valued  $n\times n$-matrix whose
$(i,j)$-entry is $0$ if the vertices  $i$ and $j$ are not connected
by an edge, it is $1$ if the arrow goes from~$i$ to~$j$, and
it is $-1$ if there is an arrow going from~$j$ to~$i$. 
The adjacency matrix of a directed graph is skew symmetric.

\begin{proposition}[\cite{DL22}]
The characteristic polynomial of the adjacency matrix
of the directed graph $\vec{g}_\alpha(C)$ is independent
of the cut point~$\alpha$. 
\end{proposition}

\begin{theorem}[\cite{DL22}]
The characteristic polynomial of the adjacency matrix of the directed 
intersection graph of a chord diagram~$C$ coincides with the 
skew characteristic polynomial $Q_{g(C)}$ of its intersection graph. 
\end{theorem}

Therefore, the skew characteristic polynomial of an intersection
graph coincides with the characteristic polynomial of the
oriented graph obtained by a certain orientation of its edges,
which belongs to a class of admissible orientations. This
fact justifies the choice of the name for the invariant. For a graph
that is not an intersection graph, there could be no such orientation.

\section{Hopf algebra of graphs}\label{s2}

Many natural graph invariants are multiplicative.
The value of such an invariant of a graph $G_1\sqcup G_2$,
which is a disjoint union of two graphs,
 $G_1$ and $G_2$, is the product of its values on~$G_1$ and~$G_2$.
In particular, all the invariants defined in  Sec.~\ref{ss4tg} are multiplicative. 
Multiplicativity means that in order to calculate the value
of the invariant on a graph it suffices to know its values
on the connected components of the graph.

In this section we show that the behaviour of the invariants 
with respect to comultiplication of graphs is of not less importance.
Multiplication and comultiplication of graphs together form
a Hopf algebra structure on the vector space spanned by simple graphs.
This Hopf algebra is a polynomial algebra in its primitive elements.
The primitive elements are linear combinations of graphs, and for many
graph invariants their values on primitive elements are much simpler
than their values of the graphs themselves. Thus, the value of a
chromatic polynomial on any primitive element is a linear polynomial,
while the chromatic polynomial of a simple graph on~$n$ vertices
is a polynomial of degree~$n$. Following~\cite{CKL20}, 
we call polynomial graph invariants whose values on primitive elements
are linear polynomials \emph{umbral invariants}. Integrability properties
of umbral invariants are discussed in Sec.~\ref{ssI}.

\subsection{Multiplication and comultiplication of graphs}

Denote by~$\cG_n$ the vector space over~$\BC$ freely spanned by
all graphs on~$n$ vertices, $n=0,1,2,\dots$, and let
$$
\cG=\cG_0\oplus\cG_1\oplus\cG_2\oplus\dots
$$
be the direct sum of these vector spaces. Note that the vector space
$\cG_0$ is one-dimensional; it is spanned by the empty graph.
Introduce a commutative multiplication $m:\cG\otimes\cG\to\cG$
on~$\cG$ by its values on the generators
$$
m:G_1\otimes G_2\mapsto G_1\sqcup G_2
$$
and extending its to linear combinations of graphs by linearity.
Obviously, this multiplication is commutative and respects the grading: 
$$
m:\cG_k\otimes \cG_n\to \cG_{k+n}
$$
for all~$k$ and~$n$. The empty graph~$e\in\cG_0$ 
is the unity of this multiplication.

The action of the \emph{comultiplication} $\mu:\cG\to\cG\otimes\cG$ 
on a graph~$G$ has the form~\footnote{Comultiplication is often
denoted by~$\Delta$; in the present paper, however, symbol~$\Delta$ 
is heavily used in the different context of delta-matroids.}
$$
\mu:G\mapsto \sum_{V(G)=U\sqcup W} G|_U\otimes G|_W,
$$
where summation in the right-hand side is carried over all
ordered partitions of the set of vertices $V(G)$ 
of~$G$ into two disjoint subsets, and $G|_U$
denotes the subgraph of~$G$ induced by a subset~$U$ of its
vertices. Comultiplication is extended to linear combinations of graphs
by linearity. With respect to this comultiplication, the vector space
$\cG$ is a 
\emph{coalgebra}. The \emph{counit} with respect to~$\mu$ 
is the linear mapping$\epsilon:\BC\to\cG_0$ taking~$1$
to~$e$.

\subsection{Primitive elements, Milnor--Moore theorem,\\ and Hopf algebra
structure}
An element~$p$ of a coalgebra with comultiplication~$\mu$ is
\emph{primitive} if
$$
\mu:p\mapsto1\otimes p+p\otimes1.
$$
In the coalgebra $\cG$, the graph~$K_1\in\cG_1$
consisting of a single vertex is primitive. No other graph
is primitive, however certain linear combinations of graphs
in which more than one graph participate are primitive.
Thus, it is easy to check that the difference  $K_2-K_1^2$
is primitive. Here and below we denote by~$K_n$ 
the complete graph on~$n$ vertices.

Any polynomial algebra $\BC[y_1,y_2,\dots]$ in either finitely or
infinitely many variables is endowed with a natural coalgebra
structure if we declare the variables $y_i$ primitive,
$\mu(y_1)=1\otimes y_i+y_i\otimes1$.
Comultiplication is extended to monomials and their linear combinations
(polynomials) as a ring homomorphism,
i.e., $\mu(\prod_{k=1}^n y_{i_k})=\prod_{k=1}^n 
(1\otimes y_{i_k}+y_{i_k}\otimes1)$.

\begin{definition}
A tuple $(H,m,\mu,e,\epsilon,S)$,  where
\begin{itemize}
    \item $H$ is a vector space;
    \item $m:H\otimes H\to H$ is a multiplication with a unit~$e:\BC\to H$;
    \item $\mu:H\to H\otimes H$ is a comultiplication with a counit $\epsilon:H\to\BC$;
    \item $S:H\to H$ is a linear mapping;
\end{itemize}
is called a \emph{Hopf algebra} if
\begin{enumerate}
    \item $m$ is a coalebra homomorphism;
    \item $\mu$ is an algebra homomorphism;
    \item $S$ satisfies the relation
    $m \circ (S \otimes {\rm Id}) \circ \mu 
    = \eta\circ\epsilon = m \circ ({\rm Id} \otimes S)
    \circ \mu$.
\end{enumerate}
\end{definition}

The mapping~$S$ in a Hopf algebra is called an \emph{antipode}.

The algebra of polynomials $\BC[y_1,y_2,\dots]$ becomes a Hopf algebra
if we define the antipode~$S$ as the mapping $S:y_i\mapsto-y_i$
for all~$i=1,2,\dots$, and extend it to whole space~$H$ as an algebra
homomorphism.

A Hopf algebra~$H$ can be \emph{graded}. In this case, it is represented
as a direct sum of finite dimensional subspaces
$$
H=H_0\oplus H_1\oplus H_2\oplus\dots,
$$
and both multiplication and comultiplication must respect the grading, i.e.,
$$
m:H_i\otimes H_j\to H_{i+j},\qquad 
\mu:H_n\to\bigoplus_{i+j=n} H_i\otimes H_j.
$$
A graded Hopf algebra is said to be \emph{connected}
if the vector space~$H_0$ is one-dimensional.

If a grading (weight) $w(y_i)\in \BN$ of each variable~$y_i$ 
in a polynomial Hopf algebra is given such that for each~$n=1,2,3,\dots$
there are finitely many variables of weight at most~$n$, then
the polynomial Hopf algebra $\BC[y_1,y_2,\dots]$ 
becomes graded. The subspace of grading~$n$ in it is spanned by
the monomials of weight~$n$, where the weight of a monomial is the sum
of the weights of the variables entering it.

The following Milnor--Moore theorem describes the structure of graded 
commutative cocom\-mutative Hopf algebras.

\begin{theorem}[\cite{MM65}]\label{tMM}
Any connected graded commutative cocommutative Hopf algebra is a
poly\-nomial Hopf algebra in its primitive elements.
\end{theorem}

This theorem means that if we choose a basis in each subspace of
primitive elements 
$P(H_n)\subset H_n\subset H$ of grading~$n$,
then $H$ is a graded polynomial Hopf algebra in the elements
of these basis if we set the weight of each element equal
to its grading.

\begin{remark}
The Milnor--Moore theorem remains true if one assumes that the Hopf
algebra is cococmmutative only. However, we are going to use it only
in the commutative case, where the statement above is sufficient.
\end{remark}

A primitive element is naturally associated to any polynomial
in a Hopf algebra of polynomials, namely, its linear part.
This mapping from the vector space of polynomials to the vector space
of their linear parts is a projection, i.e. its square coincides
with itself. Milnor--Moore theorem~\ref{tMM} implies that
the projection to the subspace of primitive elements is naturally
defined in any connected graded commutative cocommutative Hopf algebra~$H$.
We denote this projection by~$\pi$, $\pi:H\to P(H)$. 

Each homogeneous subspace $H_n\subset H$
can be represented as a direct sum $H_n=P(H_n)\oplus D(H_n)$
of the subspace of primitive elements $P(H_n)$ and the kernel
$D(H_n)$ of the projection~$\pi$. The kernel $D(H_n)$
consists of \emph{decomposable elements}, that is, of polynomials
in primitive elements of grading smaller than~$n$. 

The projection~$\pi$ can be given by the following formula.
Let $\phi,\psi:H\to H$ be linear mappings. 
Define the \emph{convolution product} 
$\phi *\psi:H\to H$ as the linear mapping 
$(\phi*\psi)(a)=m((\phi\otimes\psi)(\mu(a)))$ for all $a\in H$.
The convolution product in graded Hopf algebras 
can be used to define other operations on linear mappings,
that are represented in terms of power series.
Thus, if $\phi:H\to H$ is a grading preserving linear mapping 
of the Hopf algebra~$H$ to itself such that
$\phi(1)=1$, then the logarithm of~$\phi$ can be defined as
$$
\log~\phi=\phi_0-\frac12\phi_0*\phi_0+\frac13\phi_0*\phi_0*\phi_0-\dots,
$$
where $\phi_0:H\to H$ is the linear mapping given by $\phi_0:1\mapsto0$, 
and $\phi_0$ coincides with~$\phi$ on all homogeneous subspaces~$H_k$
of positive grading $k>0$.

\begin{theorem}[\cite{S95}]
The projection~$\pi:H\to H$ is the logarithm of the identity
mapping, $\pi=\log~{\rm Id}$.
\end{theorem}

To prove this statement, it sufficies to verify that
$\log~{\rm Id}(p)=\phi_0(p)=p$ for any primitive element~$p$, 
and $\log~{\rm Id}(p_1p_2\dots)=0$ for any nonlinear monomial in
primitive elements.

In particular, in the Hopf algebra of graphs $\cG$ the projection~$\pi$
to the subspace of primitive elements whose kernel is the subspace of
decomposable elements is given by~\cite{L97}
\begin{equation}\label{epp}
\pi(G)=G-1!\sum_{U_1\sqcup U_2=V(G)}G|_{U_1}G|_{U_1}
+2!\sum_{U_1\sqcup U_2\sqcup U_3=V(G)}G|_{U_1}G|_{U_2}G|_{U_3}
-\dots,
\end{equation}
where the summations are carried over all unordered partitions
of the set $V(G)$ of the vertices of~$G$ into $2,3,4,\dots$
nonempty disjoint subsets. For example,
$\pi(K_3)=K_3-3K_1K_2+2K_1^3$. It is easy to check that the linear combination
of graphs in the right hand side indeed is a primitive
element.

The projection $\pi:\cG\to\cG$ takes disconnected graphs to~$0$
since they are decomposable elements of the Hopf algebra~$\cG$. 
In turn, projections of connected graphs with~$n$ vertices
form a basis in the space $P(\cG_n)$ of primitive elements in grading~$n$.
 
The identity mapping ${\rm Id}:\cG\to \cG$ preserves $4$-term relations,
whence the same is true for its logarithm
$\log~{\rm Id}=\pi$. As a corollary, the projection formula works
in the Hopf algebra~$\cF$ as well.

In the Hopf algebra of chord diagrams~$\cC$,
the projection formula looks similarly, with the set $V(C)$
of chords of the chord diagram~$C$ replacing the set $V(G)$
of vertices of~$G$.

\begin{remark}
There is another way to associate to a graph a primitive elements
in the Hopf algebra of graphs, see e.g.~\cite{AM13}, which, probably,
looks simpler. Namely, let's take~$G$ to the sum
$$
G\mapsto 
\sum_{E'\subset E(G)}
(-1)^{|E(G)|-|E'|}G|_{E'},
$$
where summation is carried over all spanning subgraphs of~$G$.
However, in contrast to the projection~$\pi$
this operation is specific to the Hopf algebra of graphs and cannot
be extended by linearity to a projection to the subspace of primitives.
\end{remark}


\subsection{Graph invariants and the Hopf algebra structure}

The behaviour of many graph invariants, as well as weight systems,
is closely related to the structure of the corresponding Hopf algebra.
The chromatic polynomial demonstrates a typical example.
Let us extend the chromatic polynomial to a mapping
$\chi:\cG\to\BC[c]$ of the whole space~$\cG$
to the space of polynomials in one variable by linearity.

\begin{theorem}\label{tchrpr}
The value of a chromatic polynomial on any primitive element
of the Hopf algebra~$\cG$ is a linear polynomial,
i.e., a monomial of degree~$1$.
\end{theorem}

In fact, this property is nothing but the well known
binomiality property of chromatic polynomial, which
asserts that~$\chi$ is a Hopf algebra homomorphism:

\begin{theorem}
The chromatic polynomial of a graph~$G$ satisfies the relation
$$
\chi_G(x+y)=\sum_{U\sqcup W=V(G)}
\chi_{G|_U}(x)\chi_{G|_W}(y),
$$
where summation on the right hand side is carried over all
ordered partitions of the set~$V(G)$ of vertices of~$G$
into two disjoint nonempty subsets.
\end{theorem}

Graded homomorphisms of the Hopf algebra of graphs~$\cG$
to the Hopf algebra of polynomials $\BC[p_1,p_2,\dots]$
are studied in more detail below, see Sec.~\ref{ssI}.

Theorem~\ref{tchrpr} means, in particular, that the value of 
chromatic polynomial on the projection $\pi(G)$
of an arbitrary graph~$G$ to the subspace of primitives is a
linear polynomial. Equation~(\ref{epp}) and the fact that the free 
term of the chromatic polynomial of any (nonempty) graph is~$0$,
imply that the polynomial
$\chi_{\pi(G)}(c)$ coincides with the linear term of
$\chi_G(c)$ for any graph~$G$. 

Certain other graph invariants behave similarly.

\begin{theorem}[\cite{DL22},~\cite{K20}]
The skew characteristic polynomial of a graph on any
primitive element is a constant.

The interlace polynomial of a primitive element of order~$n$ is a polynomial of degree at most~$\left[\frac{n}2\right]$.
\end{theorem}

\subsection{Hopf algebra of graphs and integrability}\label{ssI}

In this section we describe a relation between graph invariants
and the theory of integrable systems of mathematical physics.
To a graph invariant, one can associate its averaging, namely,
the result of summation of the values of the invariant over
all graphs taken with the weight inverse to the order of the
automorphism group of graphs. If the invariant takes values in
the ring of polynomials in infinitely many variables, then the 
result is a formal power series in these variables.

For graphs on surfaces (embedded graphs) similar constructions,
as is well known, lead to solutions of integrable hierarchies.
However, for abstract graphs a result of this type has been proved
only recently~\cite{CKL20}. It asserts that the averaging of an
umbral graph invariant becomes, after an appropriate rescaling
of the variables, a solution to the Kadomtsev--Petviashvili
hierarchy.

\subsubsection{KP integrable hierarchy}

The {\it\bfseries Kadomtsev--Petviashvili integrable hierarchy {\rm(}KP hierarchy{\rm)}} is an infinite system of nonlinear partial differential equations for an unknown function $F(p_1,p_2,\dots)$ 
depending on infinitely many variables. The equations of the hierarchy
are indexed by partitions of integers $n$, $n\ge4$, into two parts
none of which equals~$1$. The first two equations, which correspond
to partitions of~$4$ and~$5$, respectively, are
\begin{eqnarray*}\frac{\partial^2F}{\partial p_2^2} &=& \frac{\partial^2F}{\partial p_1\partial p_3}-
      \frac12\Bigl(\frac{\partial^2F}{\partial p_1^2}\Bigr)^2 - \frac1{12}\frac{\partial^4F}{\partial p_1^4}\\
\frac{\partial^2F}{\partial p_2\partial p_3}&=& \frac{\partial^2F}{\partial p_1\partial p_4}-
      \frac{\partial^2F}{\partial p_1^2}
      \cdot\frac{\partial^2F}{\partial p_1\partial p_2}
			- \frac16\frac{\partial^4F}{\partial p_1^3\partial p_2}.
\end{eqnarray*}
The left hand side of each equation corresponds to a partition into 
two parts none of which equals~$1$, while the terms in the right hand side
correspond to partitions of the same number~$n$, which include 
parts equal to~$1$. For $n=6$, there are two equations corresponding
to the partitions $2+4=6$ and
$3+3=6$, and so on. Exponents of solutions to the KP hierarchy are 
called its $\tau$-functions.

\subsubsection{Umbral graph invariants}

\begin{definition}
A graph invariant with values in the ring of polynomials
$\BC[q_1,q_2,\dots]$ in infinitely many variables is called 
\emph{umbral} if its extension to a mapping
 $\cG\to\BC[q_1,q_2,\dots]$ by linearity is a graded
 Hopf algebra homomorphism; here the grading in the ring of polynomials
 is defined by the weights of the variables $w(q_i)=i$, for $i=1,2,3\dots$.
\end{definition}

The definition immediately implies that a graph invariant is umbral
iff its value on any primitive element of order~$n$
is $cq_n$ for some constant~$c$.

\begin{example}
Stanley's symmetrized chromatic polynomial is an example of umbral invariant.
Indeed, by theorem~\ref{tHawg}, the Hopf algebra of weighted graphs
modulo contraction--deletion relation is graded isomorphic to the Hopf
algebra of polynomials.
\end{example}

Other umbral invariants and their various specializations 
obtained by assigning concrete values to the variables~$q_i$
contain many well known graph invariants.

\begin{example}
In~\cite{CKL20}, Abel polynomials of graphs are introduced.
Associate to any spanning forest in a graph~$G$ the monomial
in the variables~$q_i$ equal to the product
$iq_i$ over all the trees in the forest, where~$i$ 
denotes the number of vertices in the tree. Multiplication by~$i$
here corresponds to the choice of a root in the tree,
that is, such a monomial encodes the number of rooted trees
corresponding to the chosen spanning forest. The Abel polynomial
$A_G(q_1,q_2,\dots)$ is equal to the sum of these monomials 
over all spanning forests in~$G$.

If $G$ is a graph with~$n$ vertices, then $A_G$
is a quasihomogeneous polynomial of degree~$n$, 
the coefficient of $q_n$ in it being the complexity of~$G$
times~$n$. (Recall that the
\emph{complexity} of a graph is the number of spanning trees in it.)
After the substitution $q_i=x$ for all $i=1,2,\dots,n$,  
 $A_{K_n}$ of the complete graph on~$n$ vertices, it becomes
 the conventional Abel polynomial~\cite{Ab826}
$A_n(x)=x(x+n)^{n-1}$.

\begin{theorem}[\cite{CKL20}]
The Abel polynomial is an umbral graph invariant, i.e., its 
extension to the Hopf algebra of graphs~$\cG$ 
is a graded Hopf algebra homomorphism.
\end{theorem}

\end{example}

\subsubsection{Integrability of umbral invariants}

\begin{theorem}\label{tm}
Let~$I$ be an umbral graph invariant taking values in
the ring of polynomials in infinitely many variables
 $q_1,q_2,\dots$. Define the generating function
\begin{equation}\label{eupgi}
\cI^\circ(q_1,q_2,\dots)=\sum_{G}\frac{I_G(q_1,q_2,\dots)}{|\Aut(G)|},
\end{equation}
where summation is carried over all isomorphism classes of graphs,
and $|\Aut(G)|$ denotes the order of the automorphism group of~$G$.

Denote the constants $i_n$ as the sum
$$
i_n=n!\sum_{G\text{\rm\ connected}} \frac{[q_n]I_G(q_1,q_2,\dots)}{|\Aut(G)|}
$$
of all the coefficients of the monomial~$q_n$.

Suppose the constant~$i_n$ is nonzero for all $n=1,2,\dots$.
Then after rescaling the variables~$\displaystyle q_n=\frac{2^{n(n-1)/2}(n-1)!}{i_n}\cdot p_n$
the generating function~$\cI^\circ$ becomes the following linear
combination of the one-part Schur polynomials{\rm:}
$$
S(p_1,p_2,\dots)=
1+2^0s_1(p_1)+2^1s_2(p_1,p_2)+\dots+2^{n(n-1)/2}s_n(p_1,p_2,\dots,p_n)+\dots.
$$
\end{theorem}

We underline that after the rescaling in the theorem
{\it each\/} umbral graph invariant~$I$ becomes
{\it  one and the same \/} generating function.

Summation over connected graphs in the expression for~$i_n$ 
can be replaced by summation over all graphs, since for any
disconnected graph~$G$ the coefficient of the monomial~$q_n$  in the polynomial $I_G(q_1,q_2,\dots)$ is~$0$.

Since any linear combination of one-part Schur polynomials
is a $\tau$-function of the KP hierarchy, we have

\begin{corollary}
After rescaling of the variables given in Theorem~\ref{tm},
the generating function $\cI$ becomes a
$\tau$-function for the KP hierarchy.
\end{corollary}

\section{Hopf algebra of chord diagrams\\ and weight systems associated to Lie algebras}\label{s3}

Both the space of chord diagrams and the space of graphs are equipped with natural Hopf algebra structures. We describe in this section the relationship of this structure with the weight systems constructed from Lie algebras. This is a large and important class of weight systems related to quantum knot invariants. However, finding values of weight systems in this class is computationally a rather complicated problem since it requires to make computations in a noncommutative algebra. We recall known results and describe new ones that help to overcome these difficulties and obtain explicit answers in the case of the Lie algebras $\sl(2)$ and $\gl(N)$.

\subsection{Structure of the Hopf algebra of chord diagrams}

Similarly to graphs, chord diagrams generate a Hopf algebra. An essential distinction of this Hoph algebra from the Hopf algebra of graphs is that it is well defined only in the quotient space of chord diagrams modulo the subspace spanned by the $4$-term relations. Without taking this quotient there is no way to introduce multiplication correctly (while the comultiplication is defined in a natural way).

\begin{figure}[h]
    \centerline{
\begin{tikzpicture}[baseline={([yshift=-.5ex]current bounding box.center)}]
	\draw (0,0) circle (1);
	\draw (1,0) -- (-1,0);
	\fill[black] (1,0) circle (1pt)
                 (-1,0) circle (1pt)
                 (60:1) circle (1pt)
                 (120:1) circle (1pt)
                 (240:1) circle (1pt)
                 (300:1) circle (1pt);
	\draw (60:1)  ..  controls (0,0.2) and (0,-0.2)  .. (300:1);
	\draw (120:1) ..  controls (0,0.2) and (0,-0.2)  .. (240:1);
\end{tikzpicture} $\times$
\begin{tikzpicture}[baseline={([yshift=-.5ex]current bounding box.center)}]
	\draw (0,0) circle (1);
	\draw (1,0) -- (-1,0);
	\fill[black] (1,0) circle (1pt)
                 (-1,0) circle (1pt)
                 (60:1) circle (1pt)
                 (120:1) circle (1pt)
                 (240:1) circle (1pt)
                 (300:1) circle (1pt);
	\draw (60:1)  ..  controls (0,0.2) and (0,-0.2)  .. (300:1);
	\draw (120:1) ..  controls (0,0.2) and (0,-0.2)  .. (240:1);
\end{tikzpicture} $=$
\begin{tikzpicture}[baseline={([yshift=-.5ex]current bounding box.center)}]
	\draw ([shift=( 20:1cm)]0,0) arc [start angle= 20, end angle= 340, radius=1];;
	\draw (0,1) -- (0,-1);
	\fill[black] (0,1) circle (1pt)
                 (0,-1) circle (1pt)
                 (30:1) circle (1pt)
                 (150:1) circle (1pt)
                 (210:1) circle (1pt)
                 (330:1) circle (1pt);
	\draw (30:1)  ..  controls (0.2,0) and (-0.2,0)  .. (150:1);
	\draw (210:1) ..  controls (-0.2,0) and (0.2,0)  .. (330:1);

	\draw[xshift=2.5cm] ([shift=( 200:1cm)]0,0) arc [start angle= 200, end angle=520, radius=1];
	\draw[xshift=2.5cm] (30:1) -- (210:1);
	\fill[xshift=2.5cm,black] (0,1) circle (1pt)
                 (0,-1) circle (1pt)
                 (30:1) circle (1pt)
                 (150:1) circle (1pt)
                 (210:1) circle (1pt)
                 (330:1) circle (1pt);
	\draw[xshift=2.5cm] (150:1)  ..  controls (190:0.2) and (230:0.2)  .. (0,-1);
	\draw[xshift=2.5cm] (0,1) ..  controls (50:0.2) and (10:0.2)  .. (330:1);
	\draw (20:1) -- ($(20:1)+(0.62,0)$);
	\draw (-20:1) -- ($(-20:1)+(0.62,0)$);
\end{tikzpicture} $=$
\begin{tikzpicture}[baseline={([yshift=-.5ex]current bounding box.center)}]
	\draw (0,0) circle (1);
	\draw (  0:1) -- (-90:1)
          ( 30:1) -- (-30:1)
          ( 60:1) -- (-60:1)
          ( 90:1) -- (150:1)
          (120:1) -- (210:1)
          (180:1) -- (240:1);
	\fill[black] (  0:1) circle (1pt)
                 ( 30:1) circle (1pt)
                 ( 60:1) circle (1pt)
                 ( 90:1) circle (1pt)
                 (120:1) circle (1pt)
                 (150:1) circle (1pt)
                 (180:1) circle (1pt)
                 (210:1) circle (1pt)
                 (240:1) circle (1pt)
                 (270:1) circle (1pt)
                 (300:1) circle (1pt)
                 (330:1) circle (1pt);
\end{tikzpicture}}
    \caption{Multiplication of chord diagrams}\label{fcdp}
\end{figure}
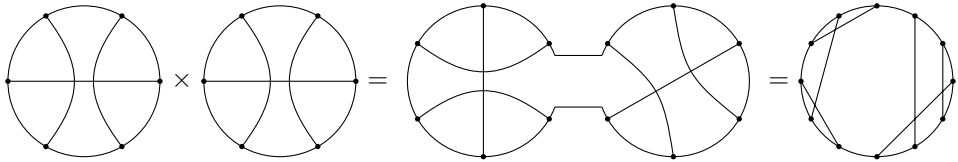

Define the product $C_1C_2$ of chord diagrams  $C_1,C_2$ as the chord diagram obtained by gluing their Wilson loops cutting them at arbitrary points distinct from the endpoint of chords, preserving the orientations of these loops, see Fig.~\ref{fcdp}. This product is well defined nodulo the $4$-term relations.

The comultiplication takes any chord diagram $C$ to the sum of tensor products of chord diagrams obtained by splitting the set of chords in~$C$ into two nonintersecting subsets
$$
\mu:C\mapsto \sum_{I\sqcup J=V(C)} C|_I\otimes C|_J.
$$
These operations turn the vector space
$$
\cC=\cC_0\oplus\cC_1\oplus\cC_2\oplus\dots,
$$
where $\cC_n$ is spanned by chord diagrams with~$n$ chords, $4$-term relations moduled out, into graded commutative cocommutative connected Hopf algebra~\cite{K93}.

The map taking a chord diagram to its intersection graph is extended by linearity to a homomorphism of the Hopf algebra of chord diagrams to the Hopf algebra of graphs modulo the $4$-term relations.

The Hopf algebra of chord diagrams is naturally isomorphic to the Hopf algebra of Jacobi diagrams. A \emph{Jacobi diagram} is an embedded $3$-regular graph with distinguished oriented loop called the Wilson loop. The grading of such a diagarm is equal to half the number of vertices (it is easy to see that the number of vertices in a $3$-regular graph is necessarily even). Denote by $\cJ_n$ the space spanned by Jakobi diagrams of grading~$n$ factorized by
\begin{itemize}
    \item skewsymmetry relations
    \item STU relations.
\end{itemize}
The skewsymmetry relation states that the change of orientation at any vertex of the diagram results in the change of the sign of the corresponding element; the STU relation is depicted on Fig.~\ref{7fig:40}. It allows one to represent any Jacobi diagram as a linear combination of chord diagrams. This correspondence establishes an isomorphism between~$\cJ$ and~$\cD$~\cite{BN95}.

Similarly to the case of chord diagrams, the product of two Jacobi diagrams is defined by gluing their Wilson loops cut at arbitrary points distinct from the vertices of the factors. After removing the Wilson loop, a Jacobi diagram splits into a number of connected components (the connected components of a chord diagram are its chords). The comultiplication takes a Jacobi diagram to the sum of tensor products of its subdiagrams obtained by splitting the set of its connected components into two subsets in all possible ways.

\begin{figure}[htbp]
\begin{picture}(100,50)(-20,10)
\thicklines
\put(30,40){\circle{30}}
\thicklines
\put(30,40){\line(0,1){15}}
\bezier{80}(30,40)(30,40)(19,29)
\bezier{80}(30,40)(30,40)(41,29)
\put(18,50){\line(1,0){24}}
\bezier{80}(95,25)(110,20)(125,25)
\put(110,23){\line(0,1){15}}
\put(110,38){\line(1,2){10}}
\put(110,38){\line(-1,2){10}}
\bezier{80}(175,25)(195,20)(215,25)
\bezier{80}(275,25)(295,20)(315,25)
\bezier{80}(187,23)(187,40)(175,55)
\bezier{80}(203,23)(203,40)(215,55)
\bezier{80}(287,23)(287,40)(315,55)
\bezier{80}(303,23)(303,40)(275,55)

\put(150,40){\makebox(0,0){$=$}}
\put(250,40){\makebox(0,0){$-$}}
\put(33,10){\makebox(0,0){(a)}}
\put(197,10){\makebox(0,0){(b)}}
\end{picture}
\caption{(a) A Jacobi diagram (b) STU relation. The halfedges at each internal vertex are oriented counterclockwise; transversal intersections of edges are not vertices
}\label{7fig:40}
\end{figure}
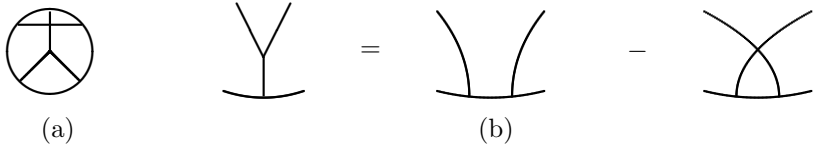

\begin{figure}[htbp]
\begin{picture}(100,40)(6,20)
\thicklines
\put(35,40){\circle{30}}
\put(105,40){\circle{30}}
\put(150,40){\circle{30}}
\put(200,40){\circle{30}}
\put(260,40){\circle{30}}
\put(310,40){\circle{30}}
\put(355,40){\circle{30}}

\put(15,40){\makebox(0,0){$\mu\Big(\,$}}
\put(68,40){\makebox(0,0){$\Big)\,=1\otimes$}}
\put(128,40){\makebox(0,0){$+$}}
\put(175,40){\makebox(0,0){$\otimes$}}
\put(228,40){\makebox(0,0){$+$}}
\put(285,40){\makebox(0,0){$\otimes$}}
\put(333,40){\makebox(0,0){$+$}}
\put(380,40){\makebox(0,0){$\otimes1$}}

\put(35,40){\line(0,1){15}}
\bezier{80}(35,40)(35,40)(24,29)
\bezier{80}(35,40)(35,40)(46,29)
\put(24,50){\line(1,0){22}}
\put(105,40){\line(0,1){15}}
\bezier{80}(105,40)(105,40)(94,29)
\bezier{80}(105,40)(105,40)(116,29)
\put(94,50){\line(1,0){22}}
\put(355,40){\line(0,1){15}}
\bezier{80}(355,40)(355,40)(344,29)
\bezier{80}(355,40)(355,40)(366,29)
\put(344,50){\line(1,0){22}}
\put(150,40){\line(0,1){15}}
\bezier{80}(150,40)(150,40)(139,29)
\bezier{80}(150,40)(150,40)(161,29)
\put(189,50){\line(1,0){22}}
\put(310,40){\line(0,1){15}}
\bezier{80}(310,40)(310,40)(299,29)
\bezier{80}(310,40)(310,40)(321,29)
\put(249,50){\line(1,0){22}}

\end{picture}
\caption{Comultiplication of Jacobi diagrams}\label{7fig:40a}
\end{figure}

In particular, connected Jacobi diagrams are primitive elements in the Hoph algebra~$\cJ$. These primitive elements span each of the homogeneous subspaces $P(\cJ_n)$ of primitive elements
(though usually they are subject to certain linear relations).
As a corollary, each of the spaces $P(\cJ_n)$, as well as the space $P(\cD_n)$ isomorphic to the former, admits a natural filtration
$$
\{0\}\subset P^{(1)}(\cJ_n)\subset P^{(2)}(\cJ_n)\subset\dots \subset P^{(n+1)}(\cJ_n)
=P(\cJ_n),
$$
where the subspace $P^{(k)}(\cJ_n)$ is spanned by connected Jacobi diagrams with at most~$k$ vertices on the Wilson loop.

The universal enveloping algebra $U\fg$ of the Lie algebra~$\fg$ admits a filtration by vector subspaces
$$
U^{(0)}\fg\subset U^{(1)}\fg\subset U^{(2)}\dots U\fg,\qquad
U^{(0)}\fg\cup U^{(1)}\fg\cup U^{(2)}\fg\cup\dots=U\fg,
$$
where $U^{(k)}\fg$ is spanned by monomials of degree at most~$k$ in the elements of~$\fg$. This filtration, in turn, induces a filtration in the center $ZU\fg$ of the universal enveloping algebra. It was shown in~\cite{CV97} that for the weight system corresponding to a given metrized Lie algebra~$\fg$, its value on a primitive element lying in the filtration term~$P^{(k)}$ belongs to $ZU^{(k)}$.

A conjecture in~\cite{ZY22} states that the projection $\pi(C)$ of a given chord diagram~$C$ with $n$ chords to the subspace of primitive elements lies in the subspace $P^{(k+1)}(\cJ_n)$, where $k$ is half the circumference of the intersection graph~$g(C)$. (The \emph{circumference} of a simple graph is the length of the longest cycle in it.) The results concerning the values of the weight systems associated to Lie algebras on the projections of chord diagrams to the subspace of primitive elements stated below confirm this conjecture.

\subsection{Universal construction of a weight system\\ associated to a Lie algebra}

To each metrized Lie algebra a weight system is associated which takes values in the center of the universal enveloping algebra of this Lie algebra. This center is typically a ring of polynomials in several generators which are the Casimir elements of the Lie algebra. We describe below the construction of this weight system and provide various algorithms for its computation leading to explicit formulas.

\subsubsection{Definition of the weight system}

The universal construction for a weight system obtained from a metrized Lie algebra looks like follows. Let $\fg$ be a Lie algebra equipped with an ${\rm ad}$-invariant scalar product (non-degenerate symmetric bilinear form) $\langle\cdot,\cdot\rangle$; here
${\rm ad}$-invariance means that the equality $\langle [a,b],c\rangle=\langle a,[b,c]\rangle$ 
holds for any three elements $a,b,c\in\fg$. Pick an arbitrary basis $e_1,\dots,e_d$ in $\fg$, $d=\dim\fg$, and denote by $e_1^*,\dots,e_d^*$ the elements of the dual basis, $\langle e_i,e_j^*\rangle=\delta_{i,j}$, $i,j=1,\dots,d$.

Now let~$C$ be a chord diagram with~$n$ chords. Choose an arbitrary point on a circle distinct from the endpoints of the chords and call it the cut point. Then, for each chord $a$, pick an index $i_a$ from the range $1,\dots,d$, and assign the basic element $e_{i_a}$ to one of the endpoints of the chord, and the element $e_{i_a}^*$ of the dual basis to the other endpoint of the chord. Multiply the resulting $2n$ elements of the Lie algebra in the order prescribed by the orientation of the Wilson loop, starting from the cut point. The resulting monomial is considered as an element of the universal enveloping algebra $U\fg$ of the given Lie algebra. Sum up the obtained monomials over all the~$d^n$ possible ways to put indices on the chords; denote the resulting sum by $w_\fg(C)\in U\fg$. For example, for the arc diagram
\begin{center}
\begin{picture}(95,40)(460,580)
\thinlines
\put(495,580){\oval(30,30)[t]}
\put(475,580){\oval(42,42)[t]}
\put(495,580){\oval(60,60)[t]}
\thicklines
\put(440,580){\vector( 1, 0){105}}
\end{picture}
\end{center}
obtained by cutting the corresponding chord diagram, the resulting element of the universal enveloping algebra is equal to
$$
\sum_{i,j,k=1}^de_ie_je_ke_i^*e_k^*e_j^*.
$$

\begin{theorem}[\cite{K93}]
\begin{enumerate}

    \item  The element $w_\fg(C)$ is independent of the choice of the basis $e_1,\dots,e_d$ in the Lie algebra.

    \item The element $w_\fg(C)$ is independent of the choice of a cut point in the chord diagram.

    \item The element $w_\fg(C)$ belongs to the center $ZU\fg$ of $U\fg$.

    \item The invariant $w_\fg$ satisfies the $4$-term relations, and thus is a weight system.

    \item The constructed weight system taking values in the commutative ring $ZU\fg$ is multiplicative, i.e.\ $w_\fg(C_1C_2)=w_\fg(C_1)w_\fg(C_2)$ for any two chord diagrams $C_1$ and $C_2$.
\end{enumerate}
\end{theorem}

\subsubsection{Jacobi diagrams and the values of the weight system on primitive elements}

The definition of the weight system associated with a Lie algebra has a convenient reformulation in terms of Jacobi diagrams. For a given metrized Lie algebra~$\fg$ consider the linear $3$-form $\omega(a,b,c)=\langle [a,b],c\rangle$ treated as an element of the third tensor power $\fg^*\otimes\fg^*\otimes\fg^*$. Denote by $\omega^*\in\fg\otimes\fg\otimes\fg$ the corresponding dual tensor obtained by identifying~$\fg$ with the dual space~$\fg^*$ using the scalar product $\langle\cdot,\cdot\rangle$. Note that the tensor $\omega^*$ is invariant with respect to cyclic permutations of the tensor factors and it changes the sign under transpositions of any two factors. Take now a Jacobi diagram. We associate the tensor $\omega^*$ with each its internal vertex labelling factors of the tensor product with the edges exiting from this vertex. To each internal edge, i.e. an edge connecting two internal vertices of the diagram, we associate the convolution of tensors corresponding to the ends of the edge by applying the scalar product. The result of such convolution is an element of the tensor power of $\fg$ whose factors are labelled by the vertices lying on the Wilson loop. Finally, we apply to this tensor the homomorphism to the universal enveloping algebra sending the tensor product of some collection of elements in~$\fg$ to their product in~$U\fg$ in the order they follow on the Wilson loop in the counterclockwise order starting from the cut point. The resulting element of~$U\fg$ is nothing but the value of the weight system~$w_{\fg}$ on the given Jacobi diagram. This definition is particularly useful for efficient computation of the values of the weight system~$w_\fg$ on primitive elements represented by connected Jacobi diagrams.

\subsection{$\sl(2)$ weight system}

The Lie algebra~$\sl(2)$ is the simplest noncommutative Lie algebra with a nondegenerate  ${\rm ad}$-invariant scalar product. In contrast to the case of more complicated Lie algebras, the $\sl(2)$ weight system admits the Chmutov-Varchenko recurrence relation, which simplifies dramatically its explicit computation. Nevertheless, even with this recursion relation the computation remains quite laborious and leads to explicit answers in a limited number of cases only. In particular, no explicit formula for the values of the~$w_{\sl(2)}$ weight system on chord diagrams all whose chords intersect one another pairwise was known until recently. We describe below the corresponding result of P.~Zakorko as well as the result of P.~Zinova and M.~Kazarian about the values of  $w_{\sl(2)}$ on chord diagrams whose intersection graph is the complete bipartite graph.

\subsubsection{Chmutov-Varchenko recursion relation}

The weight system $w_\fg$ associated with a given Lie algebra $\fg$ takes values in the commutative ring~$ZU\fg$. However, the summands of an expression for~$w_\fg(C)$ lie in the complicated noncommutative algebra~$U\fg$. Therefore, a direct application of the definition is not efficient in practice, and one needs to find methods for more efficient computation of such a weight system that would deal with commutative rings (say, rings of polynomials) at each step.

The simplest noncommutative Lie algebra is~$\sl(2)$. An efficient algorithm for computing this weight system was suggested by Chmutov and Varchenko~\cite{CV97}. The center of the universal enveloping algebra $U\sl(2)$ is the ring of polynomials in one variable, which is represented by the Casimir element. We denote this generator by~$c$. Thus, the weight system $w_{\sl(2)}$ takes values in the algebra of polynomials in~$c$.
The value $w_{\sl(2)}(C)$ on a chord diagram~$C$
with~$n$ chords is a polynomial in~$c$ of degree~$n$ with the leading coefficient~$1$ and zero free term (for $n>0$).

\begin{theorem}\label{th:Chmutov-Varchenko}
The weight system $w_{\sl(2)}$ associated to the Lie algebra $\sl(2)$ satisfies the following relations:
\begin{enumerate}
    \item $w_{\sl(2)}$ is a multiplicative weight system, $w_{\sl(2)}(C_1C_2)=w_{\sl(2)}(C_1)w_{\sl(2)}(C_2)$ for arbitrary chord diagrams $C_1$ and $C_2$. As a corollary, for the empty chord diagram the value of the weight system is equal to~$1$;
\item if the diagram $C$ is obtained from a diagram $C'$ by adding one chord having no intersections with the chords of~$C'$, then we have
$$
w_{\sl(2)}(C)=c\;w_{\sl(2)}(C').
$$
\item If the diagram $C$ is obtained from a diagram $C'$ by adding one chord intersecting exactly one chord of~$C'$, then we have
$$
w_{\sl(2)}(C)=(c-1)w_{\sl(2)}(C').
$$
\item The $6$-term relations depicted on Fig.~\ref{fCVrr} hold.
\end{enumerate}
\end{theorem}

\begin{figure}[htbp]
%
%
\def\bepi#1{\makebox[33pt]{\unitlength=18pt
            \begin{picture}(1.8,1.1)(-0.98,-0.2)  #1
            \end{picture}} }
\def\sctw#1#2#3#4{
   \bezier{25}(-0.26,0.97)(0,1.035)(0.26,0.97)
   \bezier{4}(0.26,0.97)(0.52,0.9)(0.71,0.71)
   \bezier{#1}(0.71,0.71)(0.9,0.52)(0.97,0.26)
   \bezier{4}(0.97,0.26)(1.035,0)(0.97,-0.26)
   \bezier{#2}(0.97,-0.26)(0.9,-0.52)(0.71,-0.71)
   \bezier{4}(0.71,-0.71)(0.52,-0.9)(0.26,-0.97)
   \bezier{25}(0.26,-0.97)(0,-1.035)(-0.26,-0.97)
   \bezier{4}(-0.26,-0.97)(-0.52,-0.9)(-0.71,-0.71)
   \bezier{#3}(-0.71,-0.71)(-0.9,-0.52)(-0.97,-0.26)
   \bezier{4}(-0.97,-0.26)(-1.035,0)(-0.97,0.26)
   \bezier{#4}(-0.97,0.26)(-0.9,0.52)(-0.71,0.71)
   \bezier{4}(-0.71,0.71)(-0.52,0.9)(-0.26,0.97)
}
%
\def\dslone{\bepi{\sctw{25}{25}{4}{4} \slchone \slchtwo \slchthree}}
\def\dsltwo{\bepi{\sctw{25}{25}{4}{4} \slchone \slchfour \slchthree}}
\def\dslthree{\bepi{\sctw{25}{25}{4}{4} \slchone \slchtwo \slchfive}}
\def\dslfour{\bepi{\sctw{25}{25}{4}{4} \slchone \slchfour \slchfive}}
\def\dslfive{\bepi{\sctw{25}{25}{4}{4} \slchsix \slchtth}}
\def\dslsix{\bepi{\sctw{25}{25}{4}{4} \slchseven \slcheight}}
\def\dslseven{\bepi{\sctw{4}{25}{4}{25} \slchone \slchnine \slchthree}}
\def\dsleight{\bepi{\sctw{4}{25}{4}{25} \slchone \slchten \slchthree}}
\def\dslnine{\bepi{\sctw{4}{25}{4}{25} \slchone \slchnine \slchfive}}
\def\dslten{\bepi{\sctw{4}{25}{4}{25} \slchone \slchten \slchfive}}
\def\dslel{\bepi{\sctw{4}{25}{4}{25} \slchel \put(-0.866,0.5){\circle*{0.15}}
    \put(0.866,-0.5){\circle*{0.15}} \put(0.866,-0.5){\line(-5,3){1.72}} }}
\def\dsltw{\bepi{\sctw{4}{25}{4}{25} \slchtw \slcheight}}
%
%
\def\aslone{\bepi{\sctw{25}{25}{4}{4} \slchone \slchseven \slcheight}}
\def\asltwo{\bepi{\sctw{25}{25}{4}{4} \slchone \slchseven \achone}}
\def\aslthree{\bepi{\sctw{25}{25}{4}{4} \slchone \achtwo \slcheight}}
\def\aslfour{\bepi{\sctw{25}{25}{4}{4} \slchone \achone \achtwo}}
\def\aslfive{\bepi{\sctw{25}{25}{4}{4} \slchtwo \slchthree}}
\def\asltse{\bepi{\sctw{25}{4}{25}{4} \slchone \achsix \slchtwo}}
\def\asltei{\bepi{\sctw{25}{4}{25}{4} \slchone \achsix \slchfour}}
\def\asltni{\bepi{\sctw{25}{4}{25}{4} \slchone \achfive \slchtwo}}
\def\asltte{\bepi{\sctw{25}{4}{25}{4} \slchone \achfive \slchfour}}
\def\asltto{\bepi{\sctw{25}{4}{25}{4} \achnine \slchseven}}
\def\aslttt{\bepi{\sctw{25}{4}{25}{4}\slchel\put(-0.866,-0.5){\circle*{0.15}}
   \put(0.866,0.5){\circle*{0.15}}\put(0.866,0.5){\line(-5,-3){1.72}} }}
%
%
\def\di#1{\raisebox{0pt}[20pt][23pt]{#1}}
\def\wdi#1{W\biggl( {\raisebox{0pt}[20pt][30pt]{#1}} \biggr)}
\def\spwdi#1{W\biggl( \mbox{#1} \biggr)}
\def\slrel#1#2#3#4#5#6#7{
\begin{eqnarray*}
\wdi{#1}-\wdi{#2}-\wdi{#3}+\wdi{#4} = \makebox[30pt]{}&  \\
\hfill = 2 \wdi{#5} - 2 \wdi{#6} #7 &
\end{eqnarray*} }
%
%
%
\def\slchone{\put(0,1){\circle*{0.15}}\put(0,-1){\circle*{0.15}}
             \put(0,-1){\line(0,1){2}} }
\def\slchtwo{\put(-0.174,0.985){\circle*{0.15}}
             \put(0.866,0.5){\circle*{0.15}}
             \bezier{70}(-0.174,0.985)(0.21,0.45)(0.866,0.5) }
\def\slchthree{\put(-0.174,-0.985){\circle*{0.15}}
               \put(0.866,-0.5){\circle*{0.15}}
               \bezier{70}(-0.174,-0.985)(0.21,-0.45)(0.866,-0.5) }
\def\slchfour{\put(0.174,0.985){\circle*{0.15}}
              \put(0.866,0.5){\circle*{0.15}}
              \bezier{60}(0.174,0.985)(0.36,0.48)(0.866,0.5) }
\def\slchfive{\put(0.174,-0.985){\circle*{0.15}}
              \put(0.866,-0.5){\circle*{0.15}}
              \bezier{60}(0.174,-0.985)(0.36,-0.48)(0.866,-0.5) }
\def\slchsix{\put(-0.174,0.985){\circle*{0.15}}
             \put(-0.174,-0.985){\circle*{0.15}}
             \bezier{100}(-0.174,0.985)(0.3,0)(-0.174,-0.985) }
\def\slchseven{\put(-0.174,-0.985){\circle*{0.15}}
               \put(0.866,0.5){\circle*{0.15}}
               \bezier{90}(-0.174,-0.985)(0,0)(0.866,0.5) }
\def\slcheight{\put(-0.174,0.985){\circle*{0.15}}
               \put(0.866,-0.5){\circle*{0.15}}
               \bezier{90}(-0.174,0.985)(0,0)(0.866,-0.5) }
\def\slchnine{\put(-0.174,0.985){\circle*{0.15}}
              \put(-0.866,0.5){\circle*{0.15}}
              \bezier{60}(-0.174,0.985)(-0.36,0.48)(-0.866,0.5) }
\def\slchten{\put(0.174,0.985){\circle*{0.15}}
             \put(-0.866,0.5){\circle*{0.15}}
             \bezier{90}(0.174,0.985)(-0.21,0.45)(-0.866,0.5) }
\def\slchel{\put(-0.14,0.985){\circle*{0.15}}
            \put(0.174,-0.985){\circle*{0.15}}
            \put(0.2,-0.985){\line(-1,6){0.325}} }
\def\slchtw{\put(0.174,-0.985){\circle*{0.15}}
            \put(-0.866,0.5){\circle*{0.15}}
            \bezier{90}(0.174,-0.985)(0,0)(-0.866,0.5) }
\def\slchtth{\put(0.866,0.5){\circle*{0.15}}
             \put(0.866,-0.5){\circle*{0.15}}
             \bezier{60}(0.866,0.5)(0.5,0)(0.866,-0.5) }
%
%
\def\achone{\put(0.174,0.985){\circle*{0.15}}
            \put(0.866,-0.5){\circle*{0.15}}
            \bezier{80}(0.174,0.985)(0.36,0.16)(0.866,-0.5) }
\def\achtwo{\put(0.174,-0.985){\circle*{0.15}}
            \put(0.866,0.5){\circle*{0.15}}
            \bezier{80}(0.174,-0.985)(0.36,-0.16)(0.866,0.5) }
\def\achthree{\put(-0.174,0.985){\circle*{0.15}}
            \put(-0.866,-0.5){\circle*{0.15}}
            \bezier{80}(-0.174,0.985)(-0.36,0.16)(-0.866,-0.5) }
\def\achfour{\put(-0.174,-0.985){\circle*{0.15}}
            \put(-0.866,0.5){\circle*{0.15}}
            \bezier{80}(-0.174,-0.985)(-0.36,-0.16)(-0.866,0.5) }
\def\achfive{\put(-0.174,-0.985){\circle*{0.15}}
              \put(-0.866,-0.5){\circle*{0.15}}
              \bezier{60}(-0.174,-0.985)(-0.36,-0.48)(-0.866,-0.5) }
\def\achsix{\put(0.174,-0.985){\circle*{0.15}}
            \put(-0.866,-0.5){\circle*{0.15}}
            \bezier{80}(0.174,-0.985)(-0.21,-0.45)(-0.866,-0.5) }
\def\achseven{\put(0.174,0.985){\circle*{0.15}}
             \put(0.174,-0.985){\circle*{0.15}}
             \bezier{100}(0.174,0.985)(-0.3,0)(0.174,-0.985) }
\def\acheight{\put(-0.866,0.5){\circle*{0.15}}
             \put(-0.866,-0.5){\circle*{0.15}}
             \bezier{60}(-0.866,0.5)(-0.5,0)(-0.866,-0.5) }
\def\achnine{\put(0.174,0.985){\circle*{0.15}}
             \put(-0.866,-0.5){\circle*{0.15}}
             \bezier{90}(0.174,0.985)(0,0)(-0.866,-0.5) }
%
%
\def\dslone{\bepi{\sctw{25}{25}{4}{4} \slchone \slchtwo \slchthree}}
\def\dsltwo{\bepi{\sctw{25}{25}{4}{4} \slchone \slchfour \slchthree}}
\def\dslthree{\bepi{\sctw{25}{25}{4}{4} \slchone \slchtwo \slchfive}}
\def\dslfour{\bepi{\sctw{25}{25}{4}{4} \slchone \slchfour \slchfive}}
\def\dslfive{\bepi{\sctw{25}{25}{4}{4} \slchsix \slchtth}}
\def\dslsix{\bepi{\sctw{25}{25}{4}{4} \slchseven \slcheight}}
\def\dslseven{\bepi{\sctw{4}{25}{4}{25} \slchone \slchnine \slchthree}}
\def\dsleight{\bepi{\sctw{4}{25}{4}{25} \slchone \slchten \slchthree}}
\def\dslnine{\bepi{\sctw{4}{25}{4}{25} \slchone \slchnine \slchfive}}
\def\dslten{\bepi{\sctw{4}{25}{4}{25} \slchone \slchten \slchfive}}
\def\dslel{\bepi{\sctw{4}{25}{4}{25} \slchel \put(-0.866,0.5){\circle*{0.15}}
    \put(0.866,-0.5){\circle*{0.15}} \put(0.866,-0.5){\line(-5,3){1.72}} }}
\def\dsltw{\bepi{\sctw{4}{25}{4}{25} \slchtw \slcheight}}
%
%
\def\aslone{\bepi{\sctw{25}{25}{4}{4} \slchone \slchseven \slcheight}}
\def\asltwo{\bepi{\sctw{25}{25}{4}{4} \slchone \slchseven \achone}}
\def\aslthree{\bepi{\sctw{25}{25}{4}{4} \slchone \achtwo \slcheight}}
\def\aslfour{\bepi{\sctw{25}{25}{4}{4} \slchone \achone \achtwo}}
\def\aslfive{\bepi{\sctw{25}{25}{4}{4} \slchtwo \slchthree}}
\def\asltse{\bepi{\sctw{25}{4}{25}{4} \slchone \achsix \slchtwo}}
\def\asltei{\bepi{\sctw{25}{4}{25}{4} \slchone \achsix \slchfour}}
\def\asltni{\bepi{\sctw{25}{4}{25}{4} \slchone \achfive \slchtwo}}
\def\asltte{\bepi{\sctw{25}{4}{25}{4} \slchone \achfive \slchfour}}
\def\asltto{\bepi{\sctw{25}{4}{25}{4} \achnine \slchseven}}
\def\aslttt{\bepi{\sctw{25}{4}{25}{4}\slchel\put(-0.866,-0.5){\circle*{0.15}}
   \put(0.866,0.5){\circle*{0.15}}\put(0.866,0.5){\line(-5,-3){1.72}} }}
%
%
\def\di#1{\raisebox{0pt}[20pt][23pt]{#1}}
\def\wdi#1{w_{\sl(2)}\biggl({\raisebox{0pt}[20pt][30pt]{#1}}\biggr)}
\def\spwdi#1{W\biggl( \mbox{#1} \biggr)}
\def\slrel#1#2#3#4#5#6#7{
\begin{eqnarray*}
\makebox[2pt]{} \wdi{#1}-\wdi{#2}-\wdi{#3}+\wdi{#4} 
& \\
\hfill = \,\,\, \wdi{#5} - \wdi{#6} #7 \phantom{verylongstring} &
\end{eqnarray*}
}
\slrel{\dslone}{\dsltwo}{\dslthree}{\dslfour}{\dslfive}{\dslsix};
\slrel{\aslone}{\asltwo}{\aslthree}{\aslfour}{\dslfive}{\aslfive}.
\caption{$6$-term relations for the weight system $\sl(2)$ }\label{fCVrr}
\end{figure}

As usual, it is assumed in these relations that besides the shown chords all diagrams participating in the relations contain other chords that are the same for all six diagrams. The Chmutov-Varchenko relations are sufficient to compute the value of the weight system on arbitrary chord diagram. Indeed, if the relations 1--3 are not applicable to a given chord diagram, then this diagram necessarily contains a chord to which relation 4 can be applied. All diagrams obtained by applying the relations are simpler than the original one in the following sense: the two diagrams on the right hand side have smaller number of chords than each of the diagrams on the left hand side; all the diagrams on the left hand side have equal number of chords but the last three of them have smaller number of pairs of intersecting chords than the first one. This implies that the value of the weight system is computed in a unique way by applying these relations repeatedly in finitely many steps.

The Chmutov-Varchenko relations can be also used for an axiomatic definition of the weight system $w_{\sl(2)}$. With this approach, the equalities of the theorem are taken as the definition of the weight system. Then, the assertion that the function is well defined becomes a nontrivial statement, which claims that the result of computation is independent of the order in which the Chmutov-Varchenko relations are applied. The $4$-term relations for the constructed function is a corollary of Chmutov-Varchenko relations, i.e., it is a weight system indeed.

\subsubsection{$\sl(2)$-weight system for graphs}

The following result of Chmutov and Lando is specific for the Lie algebra $\sl(2)$ and the weight system associated to it. Its analog, say, for the Lie algebra $\sl_3$ does not hold.

\begin{theorem}[\cite{CL07}]\label{conj:int}
The value of the weight system $w_{\sl(2)}$ on any chord diagram is uniquely determined by its intersection graph.
\end{theorem}

This theorem implies that the weight system $w_{\sl(2)}$ defines a function on those graphs that are intersection graphs of chord diagrams. This function vanishes on those combinations of intersection graphs that are involved into $4$-term relations. 

\begin{conjecture}[Lando]
The weight system $w_{\sl(2)}$ admits an extension to a function on graphs satisfying the $4$-term relation for graphs. This extension is unique.
\end{conjecture}\label{csl2}

Existence of such an extension would mean that $w_{\sl(2)}$ vanishes on any linear combination of intersection graphs that are corollaries of the $4$-term relations for graphs (such a linear combination is not necessarily a corollary of the $4$-term relation for chord diagrams). Uniqueness of such an extension is not specific for~$w_{\sl(2)}$. It is related to the question whether the mapping sending a chord diagram to its intersection graph is epimorphic, see Sect.~\ref{ss4tg}. Validity of the conjecture is checked using computer for graphs with at most~$8$ vertices (E.~Krasilnikov,~\cite{K21}). Note that for the extension obtained by Krasilnikov the values of $w_{\sl(2)}$ on some graphs with~$8$ vertices are not integer (all such graphs are not intersection graphs since for all intersection graphs the Chmutov--Varchenko relations guarantee that the values of the weight system on all intersection graphs are polynomials with integer coefficients).

There are also some other pieces of evidence supporting Conjecture~\ref{csl2}. For instance, the required extension is known for the top coefficient of the values of~$w_{\sl(2)}$ in the projection of the chord diagram to the subspace of primitive elements.
The value $w_{\sl(2)}(\pi(D))$ of the weight system~$w_{\sl(2)}$ on the projection of a given chord diagram~$D$
with~$2n$ chords to the subspace of primitive elements is a polynomial of degree at most~$n$.

\begin{theorem}[\cite{KLMR14},~\cite{BNV15}]
The coefficient of $c^n$ in the value $w_{\sl(2)}(\pi(D))$ of the weight system $w_{\sl(2)}$ on the projection of the chord diagram~$D$
with~$2n$ chords to the subspace of primitive elements coincides with the value $2\log~\nu(g(D))$ of doubled logarithm of the nondegeneracy of the intersection graph of the diagram (the logarithm is understood in the sense of convolution in the Hopf algebra).
\end{theorem}

Since the nondegeneracy of a graph is its $4$-invariant, we obtain a confirmation of the conjecture.

One of the possible approaches to the proof of Conjecture~\ref{csl2} consists in an attempt to construct a $4$-invariant of graphs taking values in the ring of polynomials in one variable whose values on intersection graphs coincide with those for $w_{\sl(2)}$. In order to construct such an invariant, it is important to have a large range of explicitly computed values of~$w_{\sl(2)}$ on graphs: these values might indicate which characteristics of the graph are reflected in~$w_{\sl(2)}$. In many ways, the results provided below are inspired by this problem.

In some respect, the weight system $w_{\sl(2)}$ resembles very much the chromatic polynomial of a graph: 
\begin{itemize}
\item it is multiplicative;
\item its value on a chord diagram is a polynomial in one variable whose degree is equal to the number of vertices in the intersection graph;
\item the leading coefficient of this polynomial is equal to ~$1$;
\item the signs of the coefficients alternate;
\item the absolute values of the second coefficient of the polynomial is equal to the number of edges in the graph;
\item the free term of this polynomial is equal to zero if the diagram has at least one chord;
\item adding a leaf to a chord diagram leads to the multiplication of the value of $w_{\sl(2)}$ by $c-1$.
\end{itemize}

It is worth to know, thereby, whether an analogue of a theorem by J.~Huh~\cite{H12} proved by him for the case of chromatic polynomial holds for the weight system  $w_{\sl(2)}$ as well:

\begin{problem}
Is it true that the sequence of absolute values of the coefficients of the weight system $w_{\sl(2)}$ is logarithmically convex for arbitrary chord diagram? In particular, is it true that this sequence of coefficients is unimodular, that is, it increases first and then decreases?
\end{problem}

The proof of Huh relates the chromatic polynomial to the geometry of algebraic manifolds. It would be interesting to find a similar relationship for the  $\sl(2)$ weight system.

\subsubsection{The values of the $\sl(2)$ weight system on complete graphs}

The Chmutov-Varchenko relations allow one both to compute the values of the weight system~$w_{\sl(2)}$ on specific chord diagrams and to obtain closed formulas for some infinite series of chord diagrams. In particular, consider the chord diagram~$K_n$ formed by~$n$ pairwise intersecting chords. The intersection graph for this diagram is the complete graph on~$n$ vertices. For the complete graphs the computation of one or another graph invariant causes usually no difficulty due to high symmetry of complete graphs. However, for the case of the weight system $w_{\sl(2)}$ this problem proved to be surprisingly difficult and the answer is rather unexpected.

The following statement was suggested as a conjecture by S.~Lando around 2015. It was partially proved (for the linear terms of the polynomials) in~\cite{B17}.

\begin{theorem}[\cite{Za22}]
The generating series for the values of the weight system associated to the Lie algebra~$\sl(2)$ on complete graphs admits the following infinite continuous fraction expansion
\begin{eqnarray*}
\sum_{n=0}^\infty w_{\sl(2)}(K_n)t^n&=&1+ct+c(c-1)t^2+c(c-1)(c-2)t^3
+c(c^3-6c^2+13c-7)t^4+\dots\\
&=&\frac1{1+a_0t+\frac{b_1t^2}{1+a_1t+\frac{b_2t^2}{1+a_2t+\dots}}},
\end{eqnarray*}
where
$$
a_m=-c+m\,(m+1),\qquad b_m=m^2 c-\frac{m^2(m^2-1)}{4}.
$$
\end{theorem}

It is useful to compare this decomposition to a continuous fraction with a similar decomposition for the generating function for chromatic polynomials of complete graphs:
\begin{eqnarray*}
\sum_{n=0}^\infty \chi_{K_n}(c)t^n&=&1+ct+c(c-1)t^2+c(c-1)(c-2)t^3
+c(c-1)(c-2)(c-3)t^4+\dots\\
&=&\frac1{1+a_0t+\frac{b_1t^2}{1+a_1t+\frac{b_2t^2}{1+a_2t+\dots}}},
\end{eqnarray*}
$$
a_m=-c+2m,\qquad b_m=m c-m(m-1).
$$

The proof of this theorem given in~\cite{Za22} contains also the following efficient way for computing the values of the weight system from the theorem, appearing as an intermediate step of computations. Consider the linear operator~$T$ acting on the space of polynomials in one variable~$x$ whose action is defined as follows:
$$
T(1)=x,\qquad T(x)=x^2-x,
$$
and also
$$
T(x^2f)=(2x-1)T(xf)+(2c-x-x^2)T(f)+(x-c)^2
$$
for any polynomial~$f$.

\begin{theorem}[\cite{Za22}]
The value of the weight system~$\sl(2)$ on the chord diagram with the complete intersection graph is equal to the value of the polynomial $T^n(1)$ at the point $x=c$,
$$
w_{\sl(2)}(K_n)=T^n(1)\bigm|_{x=c}.
$$
\end{theorem}

From the known value of the weight system $w_{\sl(2)}$ on complete graphs it is easy to find its values on their projections to the subspace of primitive elements.
Indeed, polynomials in complete graphs form a Hopf subalgebra in the Hopf algebra of graphs~$\cG$. 
The order of the automorphism group of the complete graph~$K_n$ on~$n$ vertices is~$n!$, 
which allows one to represent the exponential generating
function for the projections of complete graphs to the subspace of primitives as the logarithm of the exponential
generating function for complete graphs:
$$
\sum_{n=1}^\infty\pi(K_n)\frac{t^n}{n!}=
\log~\sum_{n=0}^\infty K_n\frac{t^n}{n!}.
$$
Applying $w_{\sl(2)}$ to both parts of the identity and
substituting the values $w_{\sl_2}(K_n)$ we already know,
we obtain the first terms of the expansion:
\begin{eqnarray*}
\sum_{n=1}^\infty \bar w_{\sl(2)}(K_n)\frac{t^n}{n!}&=&
c\frac{t}{1!}-c\frac{t^2}{2!}+2c\frac{t^3}{3!}+
(2c^2-7c)\frac{t^4}{4!}-(24c^2-38c)\frac{t^5}{5!}\\
&&-(16c^3-284c^2+295c)\frac{t^6}{6!}+\dots.
\end{eqnarray*}

\subsubsection{Algebra of shares}

\begin{definition}
A \emph{share} with $n$ chords is an ordered pair of
oriented intervals in which $2n$ pairwise distinct points
are given, split into~$n$ disjoint pairs, considered up to orientation preserving
diffeomorphisms of each of the two intervals.
\end{definition}

\begin{figure}
\centering
\includegraphics[scale=.6]{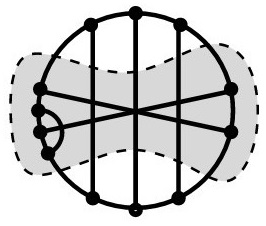}
    \caption{Join of two shares}
    \label{fig:sharejoin}
\end{figure}

For the intervals, we often take two arcs of the Wilson loop. In this case a share is a part of a chord diagram
formed by the chords whose both ends belong to a chosen 
pair of arcs, while none of the other chords has an end on these arcs.
If a pair of arcs of a chord diagram forms a share, then its complement also is a share.
In this case we call the chord diagram the 
\emph{join}   of two shares and denote it by
$S_1\#S_2$, see Fig.~\ref{fig:sharejoin}. 
In a join, the end of the first interval of the share~$S_1$
is attached to the beginning of the first interval in~$S_2$, the end of the first interval in~$S_2$ is
attached to the beginning of the second interval in~$S_1$,
the end of the second interval in~$S_1$ 
is attached to the beginning of the second interval in~$S_2$, and, finally, the end of the second interval in~$S_2$
is attached to the beginning of the first interval in~$S_1$. Join is a noncommutative operation.

To each share~$S$, one can associate a chord diagram 
by attaching the end of each of the two intervals to the
beginning of the other interval, thus uniting the 
intervals into a Wilson loop. We denote this chord diagram by~$\overline{S}$ and call it the \emph{closure}
of~$S$. We also call the intersection graph $g(\overline{S})$ the intersection graph of the share~$S$
and denote it by $g(S)$.

\begin{figure}[h]
\begin{center}
\begin{picture}(300,60)(10,40)

\thicklines
\put(40,90){\vector(-1,0){80}}
\put(-40,50){\vector(1,0){80}}
\put(-50,46){$1$}
\put(-50,86){$2$}
\put(0,50){\oval(30,30)[t]}
\put(0,50){\line(0,1){40}}
\put(-5,40){$S_1$}
\put(70,70){\circle*{4}}

\put(190,90){\vector(-1,0){80}}
\put(110,50){\vector(1,0){80}}
\put(100,46){$1$}
\put(100,86){$2$}
\put(160,90){\oval(30,30)[b]}
\put(160,90){\line(-1,-1){40}}
\put(130,90){\line(1,-1){40}}
\put(145,40){$S_2$}

\put(220,70){$=$}

\put(360,90){\vector(-1,0){100}}
\put(260,50){\vector(1,0){100}}
\put(250,46){$1$}
\put(250,86){$2$}
\put(280,50){\oval(30,30)[t]}
\put(280,50){\line(0,1){40}}
\put(340,90){\oval(30,30)[b]}
\put(340,90){\line(-1,-1){40}}
\put(300,90){\line(1,-1){40}}
\put(305,40){$S_1S_2$}

\end{picture}
\end{center}
    \caption{Product of two shares}\label{fsp}
\end{figure}
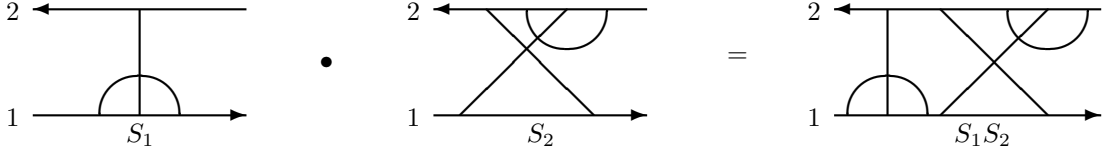

Two shares can be multiplied by attaching their intervals with coinciding numbers, see Fig.~\ref{fsp}. This
multiplication is noncommutative.

The space spanned by isomorphism classes of shares
can be endowed with relations analogous to the Chmutov--Varchenko relations in Theorem~\ref{th:Chmutov-Varchenko}. 
Namely, we can assume that all the diagrams taking part
in the relation are, in fact, shares. There is, however, some subtlety in defining relations~2 and~3 in the Theorem. Namely, suppose a share~$S$ is obtained from a share~$S'$
by adding a single chord whose both ends belong to the same
interval number $i,~i\in\{1,2\}$, which intersects at most
one of the other chords. Then the relations for the shares 
have the form
\begin{equation}\label{e1t}
S=c_i\; S',\qquad S=(c_i-1)\,S',
\end{equation}
in the cases when the given chord intersects none
of the other chords, or just one of them, respectively.

Denote by $\cS$ the quotient space of the vector space
spanned by all the shares modulo the subspace spanned
by the $6$-term relations and relations~(\ref{e1t}). 
This quotient space is endowed with the algebra structure
induced by the share multiplication operation.

Note that since the $6$-term relations are not homogeneous, the algebra~$\cS$ is filtered rather than graded.

\begin{proposition}
 The algebra $\cS$ is commutative and isomorphic to the
 algebra of polynomials in the three generators $c_1,c_2,x$ shown below:
 \begin{center}
\begin{picture}(300,60)(10,40)

\thicklines

\put(40,90){\vector(-1,0){80}}
\put(-40,50){\vector(1,0){80}}
\put(-50,46){$1$}
\put(-50,86){$2$}
\put(0,50){\oval(30,30)[t]}
\put(-5,40){$c_1$}

\put(190,90){\vector(-1,0){80}}
\put(110,50){\vector(1,0){80}}
\put(100,46){$1$}
\put(100,86){$2$}
\put(150,90){\oval(30,30)[b]}
\put(145,40){$c_2$}

\put(340,90){\vector(-1,0){80}}
\put(260,50){\vector(1,0){80}}
\put(250,46){$1$}
\put(250,86){$2$}
\put(300,90){\line(0,-1){40}}
\put(295,40){$x$}

\end{picture}
\end{center}

\end{proposition}
 
The proposition states that modulo $6$-term relations and relations~(\ref{e1t}),  each share admits a unique representation 
as a linear combination of the shares $1,x,x^2,x^3,\dots$,
where the share $x^n$ is formed by $n$ parallel chords
having ends on the two distinct intervals; the coefficients
in these linear combinations are polynomials in~$c_1$ and $c_2$.

One of the reasons for introducing the algebra $\cS$
is the following remark.

\begin{proposition}
Suppose two shares (or two linear combinations of shares)
$S_1$ and $S_2$ represent one and the same element in $\cS$. 
Then for an arbitrary share~$R$, the values of the 
$\sl_2$ weight system on the chord diagrams obtained
by joining the shares $S_1$ and $S_2$, 
respectively, with~$R$, coincide,
$$
w_{\sl(2)}(S_1\#R)=w_{\sl(2)}(S_2\#R).
$$
\end{proposition}

One can chose for the additional share~$R$, for example,
the tuple of~$n$ mutually parallel chords having ends on
distinct arcs, or the tuple of~$n$ pairwise intersecting
chords having ends on distinct arcs.

\subsubsection{Values of the $\sl(2)$-weight system on 
complete bipartite graphs}

Chord diagrams whose intersection graph is complete
bipartite form another wide class of chord diagrams,
for which the values of $w_{\sl(2)}$ are given by explicit
formulas.

Denote by $K_{m,n}$ the chord diagram formed by~$m$
parallel horizontal chords and~$n$ parallel vertical
chords. The intersection graph of such a chord diagram
is the complete bipartite graph having~$m$ vertices in one part and~$n$ vertices in the other part.
Equivalently, $K_{m,n}$ is the join of two shares formed
by, respectively, $m$ and $n$ parallel chords, each having
ends on distinct arcs. For each~$m=0,1,2,3,\dots$,
consider the generating function for the values of the
weight system $w_{\sl(2)}$ on the chord diagrams $K_{m,n}$,
$$
G_m=G_m(c,t)=\sum_{n=0}^\infty w_{\sl(2)}(K_{m,n})t^n.
$$

\begin{theorem}
The generating functions $G_{m}$ are subject to the relation
$$
\frac{G_{m}(t)-c^m}{t}=\sum_{i=0}^ms_{i,m}G_{i},\quad m=1,2,\dots,
$$
where the coefficients $s_{i,m}$ are given by the generating power series
$$
\sum_{i,m=0}^\infty s_{i,m} x^i t^m=
\frac{1}{1-x\,t}
\Biggl(c+\frac{c^2t^2-x\,t}{(1-x\,t)^2+(1+x\,t)t-2\,c\,t^2} \Biggr)
$$
\end{theorem}

Note that $G_m$ enters the right hand side of the equation
as well, with coefficient $s_{m,m}=c-\frac{m\,(m-1)}{2}$.
Moving this term to the left hand side, we can make the recurrence more explicit:
$$
G_m=\frac{1}{1-s_{m,m} t}\bigl(c^m+t\sum_{i=0}^{m-1}s_{i,m} G_i\Bigr).
$$
By induction, this implies

\begin{corollary}
The generating function $G_m$ can be represented as a finite linear combination of geometric porgressions
$\frac{1}{1-\left(c-\frac{i(i+1)}{2}\right)\,t}$, $i=0,1,\dots,m$. The coefficients of these geometric
progressions are polynomials in~$c$.
\end{corollary}
 
For small~ $m$, the recurrence yields
\begin{align*}
G_0&= \frac{1}{1-c\,t},\\
G_1&=\frac{c}{1-(c-1)\,t},\\
G_2&=\frac{c^2}{3 (1-c\,t)}
+\frac{c}{2\,(1-(c-1)\,t)}
+\frac{c\,(4 c-3)}{6\,(1-(c-3)\,t)},\\
G_3&=\frac{c^2}{6\,(1-c\,t)}
+\frac{c\,(3\,c^2-2c+2)}{5\,(1-(c-1)\,t)}
+\frac{c\,(4 c-3)}{3\,(1-(c-3)\,t)}
+\frac{c\,(c-2)(4 c-3)}{10\,(1-(c-6)\,t)}.
\end{align*}

The Proposition above admits the following generalization.
For an arbitrary share~$S$, denote by $(S,n)$ 
the chord diagram which is the join of~$S$ and the share~$x^n$ consisting of~ $n$ parallel chords.
Introduce the generating function
$$
f_S=\sum_{n=0}^\infty w_{\sl(2)}((S,n))\, t^n.
$$

\begin{corollary}
The generating function $f_S$ can be represented as a finite
linear combination of geometric progressions of the form
 $\frac{1}{1-\left(c-\frac{i(i+1)}{2}\right)\,t}$, $i=0,1,\dots,m$, where~ $m$ is the number of those chords of~$S$ whose ends belong to distinct arcs. The coefficients
 of these progressions are polynomials in~$c$.
\end{corollary}

Indeed, if a share $S$ can be presented modulo $6$-term relations and relations~(\ref{e1t}) as a linear 
combination of the basis shares $x^m$,
$$
S=\sum_{m=0}^m a_m x^m,
$$
then the series $f_S$ can be expressed as the linear combination of the series
$$
\sum_{n=0}^\infty w_{\sl(2)}(x^m,n)\, t^n
=
\sum_{n=0}^\infty w_{\sl(2)}(K_{m,n})\, t^n
=G_m,
$$
with the same coefficients, that is,
$$
f_S=\sum_{k=0}^m a_m G_m,
$$
and we can apply the Corollary above.

Similarly to the case of complete graphs, 
polynomials in complete bipartite graphs form a Hopf
subalgebra in~$\cG$. There are formulas for the values
of the weight system $w_{\sl(2)}$ on the projections
of complete bipartite graphs to the subspace of primitives,
which generalize the formula for the logarithm of the
exponential generating function for complete graphs, see
details in~\cite{F22}.

\subsubsection{Values of the $\sl(2)$-weight system
on graphs that are not intersection graphs}

Not each graph is the intersection graph of some chord diagrams, but using $4$-term relations for graphs one can
try to express it in terms of intersection graphs and extend to it the $\sl(2)$-weight system. The formulas from
the previous section allow one to obtain such an extention
not only to some specific graph, but also for infinite
sequences of graphs. For an arbitrary graph~$G$, denote by
$(G,n)$ its connected sum with~$n$ isolated vertices, i.e.,
the graph obtained by adding~$n$ vertices to~$G$ and connecting each of the added vertices to each of the vertices of~$G$. Now, take for~$G$ the graph 
$C_5=\raisebox{-0.2ex}{\includegraphics[scale=.35]{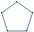}}$, the cycle on~$5$ vertices. For $n\ge1$, the graph $(C_5,n)$ is not an intersection graph.

\begin{proposition}[\cite{F20}] 
Suppose Conjecture~\ref{csl2} about extendability of the
$\sl(2)$-weight system to the space of graphs is valid.
Then the values of an extended weight system on the graphs
$(C_5,n)$ are given by the generating function
\begin{align*}
\sum_{n=0}^\infty w_{\sl(2)}((C_5,n))\,t^n\hskip-5em&\\
&=(c+5) c^2\,G_0-2 (c-1) (7 c+3)\,G_1+(5 c^2-6 c-26)\, G_2+29\,G_3-10\,G_4+G_5\\
&=\tfrac{(30 c^4-60 c^3-111 c^2+64 c+36) c}{70 (1-(c -1)t)}
+\tfrac{(c-2) (5 c^2-15 c+9) (4 c-3) c}{45 (1-(c -6) t)}
+\tfrac{(c-6) (c-2) (4 c-15) (4 c-3) c}{126 (1-(c-15) t)},
\end{align*}
where $G_m=\sum_{n=0}^\infty w_{\sl(2)}(K_{m,n})\,t^n$ are the generating functions for the values of the
$\sl(2)$-weight system on complete  bipartite graphs.
\end{proposition}

\begin{figure}
\centering
\begin{picture}(350,80)
\put(0,10){\includegraphics[scale = 0.8]{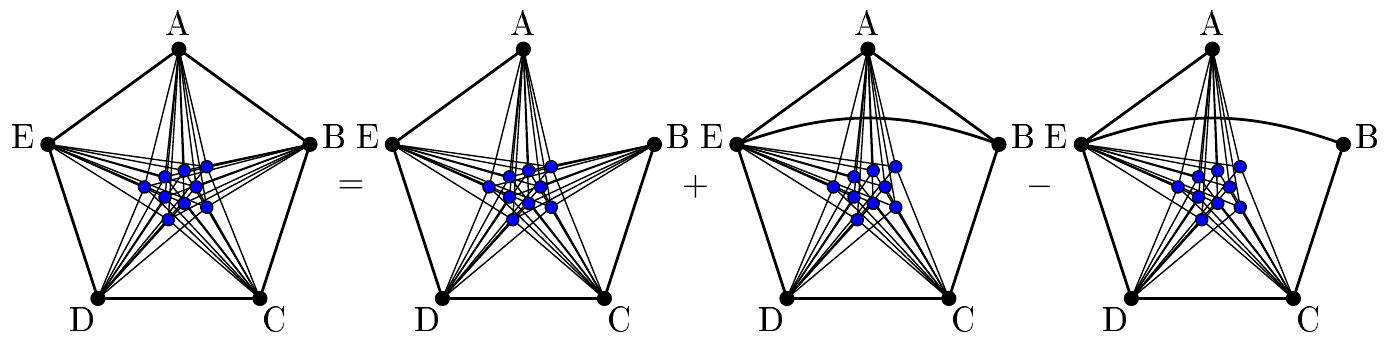}}
\put(30,0){$(C_5,n)$}
\put(110,0){$G_{1,n}$}
\put(195,0){$G_{2,n}$}
\put(270,0){$G_{3,n}$}
\end{picture}
\caption{A $4$-term relation for the graphs $(C_5,n)$.
All the graphs in the right hand side are intersection graphs} \label{pic:4term5wheel}
\end{figure}
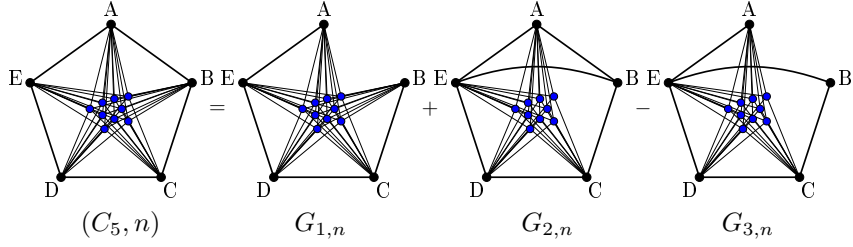

The proof is achieved by applying 
to the graphs $(C_5,n)$ the $4$-term relation  shown in Fig.~\ref{pic:4term5wheel},
$$
(C_5,n)=G_{1,n}-G_{2,n}+G_{3,n}.
$$
Each of the graphs in the right hand side is an intersection graph. Moreover, it is the intersection graph
of a chord diagram of the form $(S,n)$, for some share $S$,
$$
G_{1,n}=g\Biggl(\raisebox{-3ex}{\includegraphics[scale=0.3]{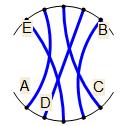}},n\Biggr),
\qquad G_{2,n}=g\Biggl(\raisebox{-3ex}{\includegraphics[scale=0.3]{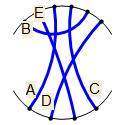}},n\Biggr),
\qquad G_{3,n}=g\Biggl(\raisebox{-3ex}{\includegraphics[scale=0.3]{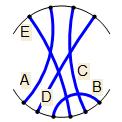}},n\Biggr).
$$
By applying the argument from the previous section, we
can express each of the shares obtained in terms of the
basic shares $x^m$, $m=0,1,\dots,5$, and express in this
way the generating function for the graphs
$(C_5,n)$ in terms of the corresponding generating functions for complete bipartite graphs.

\subsection{$\gl(N)$-weight system}

Until recently, no recursion similar to the Chmutov--Varchenko one has been known for Lie algebras
other than $\sl(2)$, and essentially the only known way of computing the values of the weight system associated
to this or that Lie algebra was the straightforward application of the definition, which is extremely laborious. Here we explain the results of the recent paper~\cite{ZY22} devoted to computation of the weight
systems associated to the Lie algebras $\gl(N)$.

For the basis in $\gl(N)$, one can choose the matrix units
$E_{ij}$, which have~$1$ on the intersection of the
$i$~th row and $j$~th column, all whose other entries are zeroes. The commutator of two matrix units has the form $[E_{ij},E_{kl}]=\delta_{jk}E_{il}-\delta_{il}E_{jk}$.
The standard invariant scalar product is given by
$\langle x,y\rangle=\mathop{\rm Tr}(xy)$,
and the dual basis with respect to this scalar product is
$E_{ij}^*=E_{ji}$. The center $ZU(\gl(N))$ is
generated freely by the  \emph{Casimir elements} $C_1,C_2,\dots,C_N$, where  $C_{1}=\sum_{i=1}^NE_{ii}$, $C_{2}=\sum_{i,j=1}^NE_{ij}E_{ji}$,
and, more generally,
$$
C_{k}=\sum_{i_1,i_2,\dots,i_k=1}^NE_{i_1i_2}E_{i_2i_3}\dots E_{i_k i_1}.
$$
Casimir elements with the numbers greater than~$N$ are
defined similarly. They also belong to the center, but they can be expressed as polynomials in the Casimir elements whose numbers are at most~$N$.
Thus, the value of the weight system associated to the Lie algebra $\gl(N)$ is a polynomial in the Casimir elements.
For example, it is easy to see that for the chord diagram~$K_1$, which consists of the single chord, we have
$w_{\gl(N)}(K_1)=C_2$.

\begin{theorem}[\cite{ZY22}]\label{th:wgl}
There is a universal weight system, which we denote by
$w_{\gl}$, taking values in the ring of polynomials in
infinitely many variables $N,C_1,C_2,\dots$, 
whose value for a given positive integer~$N$ 
coincides with the value of the $\gl(N)$ weight system. 
\end{theorem}

It is easy to see that for each chord diagram~$C$ with~$n$ chords the value of the weight system $w_{\gl(N)}(C)$, for $N\ge n$, is a polynomial in $C_1,\dots,C_{n}$. 
The theorem asserts that the coefficients of this polynomial are polynomials in~$N$. Moreover, the universal formula for $w_{\gl(N)}(C)$ we obtain remains valid for  $N< n$ as well.

\begin{example}
For the chord diagram $K_2$, which is formed by two intersecting chords, we have 
$$
w_{\gl(N)}(K_2)=\sum_{i,j,k,l=1}^NE_{ij}E_{kl}E_{ji}E_{lk}
=C_2^2+C_1^2-N\, C_2=w_{\gl}(K_2).
$$
\end{example}

Below, we define the weight system $w_{\gl}$ for arbitrary permutations and give a recurrence relation for its computation. 

\subsubsection{$\gl(N)$-weight system on permutation and recurrence relations}

Let $m$ be a positive integer, and let $S_m$ be the group of all permutations of the elements $\{1,2,\dots,m\}$;
for an arbitrary $\sigma\in S_m$, we set
$$
w_{\gl(N)}(\sigma)=\sum_{i_1,\dots,i_m=1}^NE_{i_1i_{\sigma(1)}}E_{i_2i_{\sigma(i_2)}}\dots E_{i_mi_{\sigma(m)}}\in U(\gl(N)).
$$

The element thus constructed belongs, in fact, to the center
$ZU(\gl(N))$. In addition, it is invariant under conjugation by the cyclic permutation,
$$
w_{\gl(N)}(\sigma)=\sum_{i_1,\dots,i_m=1}^N E_{i_2i_{\sigma(i_2)}}\dots E_{i_mi_{\sigma(m)}}E_{i_1i_{\sigma(1)}}.
$$

For example, the Casimir element $C_m$ corresponds to the cyclic permutation
$1\mapsto2\mapsto\dots\mapsto m\mapsto 1$.
On the other side, any arc diagram with~$n$ arcs can be considered as an involution without fixed points on the
set of $m=2n$ elements. The value of $w_{\gl(N)}$ on a
permutation given by this involution coincides with
the value of the $\gl(N)$-weight system on the corresponding chord diagram. For example, for the chord
diagram~$K_n$ we have $w_{\gl(N)}(K_n)=w_{\gl(N)}((1~n{+}1)(2~n{+}2)\dots(n~2n))$.

Any permutation can be represented by a directed graph.
The~$m$ vertices of this graph correspond to the permuted elements. They are situated along the oriented circle and are ordered cyclically in the direction of the circle. 
Directed edges of the graph show the action of the permutation (so that at each vertex there is one incoming and one outgoing edge). The directed graph 
$G(\sigma)$ consists of these $m$ vertices and $m$ edges, for example,
	\[
	G(\left(1\ n+1)(2\ n+2)\cdots(n\ 2n)\right)=\begin{tikzpicture}[baseline={([yshift=-.5ex]current bounding box.center)},decoration={markings, mark= at position .55 with {\arrow{stealth}}}]
		\draw[->,thick] (-2,0)--(2,0);
		\draw[blue,postaction={decorate}] (-1.8,0) ..  controls (-1,.5) ..(.2,0);
		\draw[blue,postaction={decorate}] (-1.4,0) ..  controls (-.6,.5) ..(.6,0);
		\draw[blue,postaction={decorate}] (-.2,0) ..  controls (1,.5) ..(1.8,0);
		\draw[blue,postaction={decorate}] (.2,0) ..  controls (-1,-.5) ..(-1.8,0);
		\draw[blue,postaction={decorate}] (.6,0) ..  controls (-.6,-.5) ..(-1.4,0);
		\draw[blue,postaction={decorate}] (1.8,0) ..  controls (1,-.5) ..(-.2,0);
		\fill[black] (-1.8,0) circle (1pt) node[below] {\tiny 1};
		\fill[black] (-1.4,0) circle (1pt) node[below] {\tiny 2};
		\fill[black] (-0.2,0) circle (1pt) node[below] {\tiny n};
		\fill[black] ( .15,0) circle (1pt) node[below] {\tiny n+1};
		\fill[black] ( .6,0) circle (1pt)  node[below]{\tiny n+2};
		\fill[black] ( 1.8,0) circle (1pt) node[below] {\tiny 2n};
		\node[below] at (-.8,0) {$\cdots$};
		\node[below] at (1.2,0) {$\cdots$};
\end{tikzpicture}
	\]
(the circle containing the vertices of the graph is shown as the horizontal line).

\begin{theorem}[\cite{ZY22}]
The value of the weight system $w_{\gl(N)}$ on a permutation possesses the following properties:
\begin{itemize}
    \item the weight system $w_{\gl(N)}$ 
is multiplicative with respect to the connected sum
(concatenation) of permutations. It follows, therefore, that for the empty graph (having zero vertices)
the value of $w_{\gl(N)}$ is~$1$.
    \item For the cyclic permutation whose cyclic order is consistent with the permutation the value of $w_{\gl(N)}$ coincides with the corresponding Casimir element, $w_{\gl(N)}(1\mapsto2\mapsto\dots\mapsto m\mapsto 1)=C_m$.
    \item  For an arbitrary permutation $\sigma\in S_m$,
    and for any two neighboring elements $k$, $k+1$ 
    in the set of vertices $\{1,2,\dots,m\}$, 
    the values of the invariant  $w_{\gl_N}$ is subject to the identity
\begin{multline*}
	w_{\gl(N)}\left(\begin{tikzpicture}[baseline={([yshift=-.5ex]current bounding box.center)},decoration={markings, mark= at position .55 with {\arrow{stealth}}}]
		\draw[->,thick] (-1,0) --  (1,0);
		\fill[black] (-.3,0) circle (1pt) node[below] {\tiny k};
		\fill[black] (.3,0) circle (1pt) node[below] {\tiny k+1};
		\draw (-.5,.8) node[left] {a};
		\draw (-.5,-.8) node[left] {b};
		\draw (.5,.8) node[right] {c};
		\draw (.5,-.8) node[right] {d};
		\draw[blue,postaction={decorate}] (-.5,.8) -- (.3,0);
		\draw[blue,postaction={decorate}] (-.3,0) -- (.5,.8);
		\draw[blue,postaction={decorate}] (-.5,-.8) -- (-.3,0);
		\draw[blue,postaction={decorate}] (.3,0) -- (.5,-.8);
	\end{tikzpicture}\right)-
		w_{\gl(N)}\left(\begin{tikzpicture}[baseline={([yshift=-.5ex]current bounding box.center)},decoration={markings, mark= at position .55 with {\arrow{stealth}}}]
		\draw[->,thick] (-1,0) --  (1,0);
		\fill[black] (.3,0) circle (1pt) node[below] {\tiny k+1};
		\fill[black] (-.3,0) circle (1pt) node[below] {\tiny k};
		\draw (-.5,.8) node[left] {a};
		\draw (-.5,-.8) node[left] {b};
		\draw (.5,.8) node[right] {c};
		\draw (.5,-.8) node[right] {d};
		\draw[blue,postaction={decorate}] (-.5,.8) -- (-.3,0);
		\draw[blue,postaction={decorate}] (.3,0) -- (.5,.8);
		\draw[blue,postaction={decorate}] (-.5,-.8) -- (.3,0);
		\draw[blue,postaction={decorate}] (-.3,0) -- (.5,-.8);
	\end{tikzpicture}\right)\\=
		w_{\gl(N)}\left(\begin{tikzpicture}[baseline={([yshift=-.5ex]current bounding box.center)},decoration={markings, mark= at position .55 with {\arrow{stealth}}}]
		\draw[->,thick] (-1,0)  -- (1,0);
		\fill[black] (0,0) circle (1pt) node[above] {\tiny k'};
		\draw (-.5,.8) node[left] {a};
		\draw (-.5,-.8) node[left] {b};
		\draw (.5,.8) node[right] {c};
		\draw (.5,-.8) node[right] {d};
		\draw[blue,postaction={decorate}] (-.5,.8) ..controls (0,.4) .. (.5,.8);
		\draw[blue,postaction={decorate}] (-.5,-.8) -- (0,0);
		\draw[blue,postaction={decorate}] (0,0) -- (.5,-.8);
	\end{tikzpicture}\right)-
		w_{\gl(N)}\left(\begin{tikzpicture}[baseline={([yshift=-.5ex]current bounding box.center)},decoration={markings, mark= at position .55 with {\arrow{stealth}}}]
		\draw[->,thick] (-1,0) --  (1,0);
		\fill[black] (0,0) circle (1pt) node[below] {\tiny k'};
		\draw (-.5,.8) node[left] {a};
		\draw (-.5,-.8) node[left] {b};
		\draw (.5,.8) node[right] {c};
		\draw (.5,-.8) node[right] {d};
		\draw[blue,postaction={decorate}] (-.5,-.8) ..controls (0,-.4) .. (.5,-.8);
		\draw[blue,postaction={decorate}] (-.5,.8) -- (0,0);
		\draw[blue,postaction={decorate}] (0,0) -- (.5,.8);
	\end{tikzpicture}\right)
\end{multline*}

The graphs on the left hand side of the identity show two neighboring vertices and the edges incident to them.
In the graphs on the right hand side, these two vertices are replaced with a single one. All the other vertices and the (semi)edges are the same for all the four graphs
participating in the identity.

In the exceptional case $\sigma(k+1)=k$, the identity acquires the form 
\begin{multline*}
		w_{\gl(N)}\left(\begin{tikzpicture}[baseline={([yshift=-.5ex]current bounding box.center)},decoration={markings, mark= at position .55 with {\arrow{stealth}}}]
		\draw[->,thick] (-1,0) --  (1,0);
		\fill[black] (-.3,0) circle (1pt) node[below] {\tiny k};
		\fill[black] (.3,0) circle (1pt) node[below] {\tiny k+1};
		\draw (-.5,.8) node[left] {a};
		\draw (.5,.8) node[right] {b};
		\draw[blue,postaction={decorate}] (-.5,.8) -- (.3,0);
		\draw[blue,postaction={decorate}] (-.3,0) -- (.5,.8);
		\draw[blue,postaction={decorate}] (.3,0) ..controls(0,-.3).. (-.3,0);
	\end{tikzpicture}\right)-
		w_{\gl(N)}\left(\begin{tikzpicture}[baseline={([yshift=-.5ex]current bounding box.center)},decoration={markings, mark= at position .55 with {\arrow{stealth}}}]
		\draw[->,thick] (-1,0) --  (1,0);
		\fill[black] (.3,0) circle (1pt) node[below] {\tiny k+1};
		\fill[black] (-.3,0) circle (1pt) node[below] {\tiny k};
		\draw (-.5,.8) node[left] {a};
		\draw (.5,.8) node[right] {b};
		\draw[blue,postaction={decorate}] (-.5,.8) -- (-.3,0);
		\draw[blue,postaction={decorate}] (.3,0) -- (.5,.8);
		\draw[blue,postaction={decorate}] (-.3,0) ..controls(0,-.3).. (.3,0);
	\end{tikzpicture}\right)\\
	=C_1\times
		w_{\gl(N)}\left(\begin{tikzpicture}[baseline={([yshift=-.5ex]current bounding box.center)},decoration={markings, mark= at position .55 with {\arrow{stealth}}}]
		\draw[->,thick] (-1,0)  -- (1,0);
		\draw (-.5,.8) node[left] {a};
		\draw (.5,.8) node[right] {b};
		\draw[blue,postaction={decorate}] (-.5,.8) ..controls (0,.4) .. (.5,.8);
	\end{tikzpicture}\right)
	-N\times
		w_{\gl(N)}\left(\begin{tikzpicture}[baseline={([yshift=-.5ex]current bounding box.center)},decoration={markings, mark= at position .55 with {\arrow{stealth}}}]
		\draw[->,thick] (-1,0) --  (1,0);
		\fill[black] (0,0) circle (1pt) node[above] {\tiny k'};
		\draw (-.5,.8) node[left] {a};
		\draw (.5,.8) node[right] {b};
		\draw[blue,postaction={decorate}] (-.5,.8) -- (0,0);
		\draw[blue,postaction={decorate}] (0,0) -- (.5,.8);
	\end{tikzpicture}\right)
\end{multline*}

\end{itemize}
\end{theorem}

Applying the identity in the theorem each graph can be reduced to a monomial in the generators
$C_k$ (that is, concatenation of independent cycles) modulo
graphs with smaller number of vertices. This leads to an inductive computation of the values of the invariant $w_{\gl(N)}$. 

\begin{corollary}
The value of $w_{\gl(N)}$ on an arbitrary permutation belongs to the center $ZU\gl(N)$ of the universal enveloping algebra, is a polynomial in $N,C_1,C_2,\dots$,
and this polynomial is universal (one and the same for all
the $\gl(N)$).
\end{corollary}

We denote the mapping taking a permutation to this universal polynomial by $w_{\gl}$. Theorem~\ref{th:wgl}
specializes this Corollary to the case when the permutation is an involution without fixed points.

Below, we give a table of values of $w_{\gl}$ on chord diagrams $K_n$, for few first values of $n$.
\begin{align*}
w_{\gl}(K_2)&=C_1^2+C_2^2-N\,C_2,\\
w_{\gl}(K_3)&=C_2^3+3 C_1^2 C_2+2 C_2 N^2+(-2 C_1^2-3 C_2^2) N,\\
w_{\gl}(K_4)&=
  3 C_1^4-4 C_1^3+6 C_2^2 C_1^2+2 C_1^2-8 C_3 C_1+C_2^4+6 C_2^2+(-6 C_2^3-14 C_1^2 C_2+6 C_1 C_2-2 C_2+2 C_4) N
\\&\qquad 
 +(6 C_1^2+11 C_2^2-2 C_3) N^2-6 C_2 N^3,\\
w_{\gl}(K_5)&=
  C_2^5+10 C_1^2 C_2^3+30 C_2^3+15 C_1^4 C_2-20 C_1^3 C_2+10 C_1^2 C_2-40 C_1 C_3 C_2+(-20 C_1^4+48 C_1^3
\\&\qquad 
  -50 C_2^2 C_1^2-32 C_1^2+30 C_2^2 C_1+96 C_3 C_1-10 C_2^4-82 C_2^2+10 C_2 C_4) N+(35 C_2^3+70 C_1^2 C_2
\\&\qquad 
  -72 C_1 C_2-10 C_3 C_2+32 C_2-24 C_4) N^2+(-24 C_1^2-50 C_2^2+24 C_3) N^3+24 C_2 N^4
\end{align*}

Similarly to the case of $\sl(2)$, the values of the weight system $w_{\gl(N)}$ on projections of complete graphs
to the subspace of primitives can be computed by taking
the logarithm of the corresponding exponential generating function.

\subsubsection{Computing the $\gl(N)$-weight system\\ by means of the Harish-Chandra isomorphism}

In this section we describe one more way, suggested in~\cite{ZY22}, to compute the $\gl(N)$ weight system.
The algebra $U(\gl(N))$ admits the decomposition into the direct sum
\begin{equation}\label{HCd}
	U(\gl(N))=(\fn_-U(\gl(N))+U(\gl(N))\fn_+)\oplus U(\fh),
\end{equation}
where $\fn_-$, $\fh$ and $\fn_+$ are, respectively, the subalgebras of low-triangular, diagonal, and upper-triangular matrices in $\gl(N)$.

\begin{definition}
The \emph{Harish-Chandra projection} is the linear projection to the second summand in~\eqref{HCd},
$$
\phi:U(\GL_N)\to U(\fh)=\BC[E_{11},\cdots,E_{NN}].
$$
\end{definition}

\begin{theorem}[Harish-Chandra isomorphism, \cite{O91}]
The restriction of the Harish-Chandra projection to the center $ZU(\GL_N)$ is an injective homomorphism of 
commutative algebras, and its isomorphic image consists of
symmetric functions in the shifted matrix units $x_i=E_{ii}+N-i$, $i=1,\dots,N$.
\end{theorem}

The isomorphism of the theorem allows one to identify the codomain
$\BC[C_1,C_2,\dots,C_N]$ of the weight system $w_{\gl(N)}$
with the ring of symmetric functions in the generators
$x_1,\dots,x_N$. The generators $C_k$ under this isomorphism admit the following explicit expression:
$$
1-N\,u-\sum_{k=1}^\infty\phi(C_k)\,u^{k+1}=\prod_{i=1}^N\frac{1-(x_i+1)\,u}{1-x_iu}.
$$

Computation of the values of the $\gl(N)$ weight system on a given chord diagram, by means of the Harish-Chandra isomorphism
consists in applying the projection $\phi$ to each monomial entering the definition of the weight system one by one and obtain its contribution to the final value. The total
contribution is obtained by summing up elements of a commutative ring. As a result, one can considerably decrease the required amount of random-access memory.
Nevertheless, time resources grows rather rapidly as~$N$ grows, since the number of monomials entering the definition, which is $N^{2n}$ for a chord diagram with~$n$ chords, also grows. Ramark also that the Harish-Chandra 
projection can be applied to computing the
$w_{\gl(N)}$ weight system for a specific~$N$, and it does not lead to a universal (polynomial) dependence on~$N$ for its values.

\subsection{$\sl(N)$ weight system}
The Lie algebra $\gl(N)$ is not simple, it can be represented as the direct sum of the simple Lie algebra
$\sl(N)$ and a one-dimensional Abelian Lie algebra.
Therefore, the center of the universal enveloping algebra
$\gl(N)$ is the tensor product of the centers of the universal enveloping algebras of $\sl(N)$ and $\BC$, 
whence it could be identified with the ring of polynomials
in $C_1$ with coefficients in $ZU(\sl(N))$. 
Therefore, the values of the weight system  $w_{\sl(N)}$
can be obtained from that of $w_{\gl(N)}$ by setting $C_1=0$ and $C_k=\tilde C_k$, $k\ge2$, where $\tilde C_k$
is the projection of the corresponding Casimir element in
$ZU(\gl(N))$ to $ZU(\sl(N))$. The result of this substitution is a polynomial in $\tilde C_2,\tilde C_3,\dots$. More explicitly, set $\tilde E_{i,j}=E_{i,j}-\delta_{i,j}N^{-1}C_1\in\sl(N)\subset\gl(N)$. Then
$$
\tilde C_{k}=\sum_{i_1,i_2,\dots,i_k=1}^N\tilde E_{i_1i_2}\tilde E_{i_2i_3}\dots \tilde E_{i_k i_1}=\sum_{i=0}^k\binom{k}{i}(-1)^iC_{k-i}\left(\tfrac{C_1}{N}\right)^i,
$$
where we set $C_0=N$ and
$$
ZU(\sl(N))=\BC[\tilde C_2,\dots,\tilde C_N]\subset
ZU(\gl(N))=\BC[C_1,\dots,C_N].
$$

Alternatively, the weight system $w_{\sl(N)}$ can be computed by extending it to permutations by means of recurrence relations similar to that for the case
$\gl(N)$: just replace there $w_{\gl(N)}$ with $w_{\sl(N)}$, $C_k$ with $\tilde C_k$ and set $\tilde C_1=0$. Similarly to the $w_{\gl}$ weight system, the result
is a universal polynomial in $N,\tilde C_2,\tilde C_3,\dots$, which we denote by $w_{\sl}$. Since the number
of generators is smaller, computation of 
$w_{\sl(N)}$ is more efficient, and the answer is more compact.
For example, for the chord diagrams $K_n$ we obtain
\begin{align*}
w_{\sl}(K_2)&=\tilde C_2^2-\tilde C_2 N,\\
w_{\sl}(K_3)&=\tilde C_2^3-3 \tilde C_2^2 N+2 \tilde C_2 N^2,\\
w_{\sl}(K_4)&=\tilde C_2^4+6 \tilde C_2^2-2 (3 \tilde C_2^3+\tilde C_2-\tilde C_4) N+(11 \tilde C_2^2-2 \tilde C_3) N^2-6 \tilde C_2 N^3,\\
w_{\sl}(K_5)&=
  \tilde C_2^5+30 \tilde C_2^3-2 (5 \tilde C_2^4+41 \tilde C_2^2-5 \tilde C_4 \tilde C_2) N+(35 \tilde C_2^3-10 \tilde C_3 \tilde C_2+32 \tilde C_2-24 \tilde C_4) N^2
\\&\qquad
  -2 (25 \tilde C_2^2-12 \tilde C_3) N^3+24 \tilde C_2 N^4,\\
w_{\sl}(K_6)&=
  \tilde C_2^6+90 \tilde C_2^4+264 \tilde C_2^2-240 \tilde C_4 \tilde C_2+160 \tilde C_3^2+(-15 \tilde C_2^5-552 \tilde C_2^3+30 \tilde C_4 \tilde C_2^2+64 \tilde C_3 \tilde C_2-72 \tilde C_2
\\&\qquad
  +88 \tilde C_4-16 \tilde C_6) N+(85 \tilde C_2^4-30 \tilde C_3 \tilde C_2^2+1014 \tilde C_2^2-174 \tilde C_4 \tilde C_2-88 \tilde C_3+32 \tilde C_5) N^2
\\&\qquad
  +(-225 \tilde C_2^3+174 \tilde C_3 \tilde C_2-416 \tilde C_2+224 \tilde C_4) N^3+2 (137 \tilde C_2^2-120 \tilde C_3) N^4-120 \tilde C_2 N^5,\\
w_{\sl}(K_7)&=
  \tilde C_2^7+210 \tilde C_2^5+3192 \tilde C_2^3-1680 \tilde C_4 \tilde C_2^2+1120 \tilde C_3^2 \tilde C_2+(-21 \tilde C_2^6-2212 \tilde C_2^4+70 \tilde C_4 \tilde C_2^3+448 \tilde C_3 \tilde C_2^2
\\&\qquad
  -10680 \tilde C_2^2+7432 \tilde C_4 \tilde C_2-112 \tilde C_6 \tilde C_2-4096 \tilde C_3^2) N+(175 \tilde C_2^5-70 \tilde C_3 \tilde C_2^3+8358 \tilde C_2^3-714 \tilde C_4 \tilde C_2^2
\\&\qquad
  -2792 \tilde C_3 \tilde C_2+224 \tilde C_5 \tilde C_2+3456 \tilde C_2-3392 \tilde C_4+544 \tilde C_6) N^2+(-735 \tilde C_2^4+714 \tilde C_3 \tilde C_2^2-12892 \tilde C_2^2
\\&\qquad
  +2212 \tilde C_4 \tilde C_2+3392 \tilde C_3-1088 \tilde C_5) N^3+4 (406 \tilde C_2^3-581 \tilde C_3 \tilde C_2+1316 \tilde C_2-464 \tilde C_4) N^4
\\&\qquad
  -12 (147 \tilde C_2^2-200 \tilde C_3) N^5+720 \tilde C_2 N^6.
\end{align*}

The values of the weight system $w_{\gl}$ can be reconstructed from that of $w_{\sl}$ by using the Hopf
algebra structure of the space of chord diagrams. Namely,
for any \emph{primitive} element $P$ of degree greater than $1$ in the algebra of chord diagrams we have $w_{\gl(N)}(P)=w_{\sl(N)}(P)\in ZU(\sl(n))$, i.e., $w_{\gl(N)}(P)$ can be expressed as a polynomial in $\tilde C_2,\tilde C_3,\dots$. For $P=K_1$, we have 
$$
w_{\gl(N)}(K_1)=C_2=\tilde C_2+\frac{C_1^2}{N}=w_{\sl(N)}(K_1)+\frac{C_1^2}{N}.
$$ 
This allows one to reconstruct the values of $w_{\gl}$ 
from the known values of $w_{\sl}$ by presenting a given chord diagram as a polynomial in primitive elements.
More explicitly, for a given chord diagram~$C$ we have
$$
w_{\gl}(C)=\sum_{I\sqcup J=V(C)} \Bigl(\tfrac{C_1^2}{N}\Bigr)^{|I|} w_{\sl}(C|_J).
$$
In particular, exponential generating functions for the
values of $w_{\gl}$ and $w_{\sl}$ on the chord diagrams $K_n$ are subject to the relation
$$
1+\sum_{n=1}^\infty w_{\gl}(K_n)\frac{x^n}{n!}=e^{\frac{C_1^2}{N}}\Biggl(1+\sum_{n=1}^\infty w_{\sl}(K_n)\frac{x^n}{n!}\Biggr).
$$

\subsection{The weight system $\gl(1|1)$}

In addition to metrized Lie algebras, weight systems can be constructed from metrized Lie superalgebras. Such a weight system takes values on the center of the universal enveloping algebra of the corresponding Lie superalgebra.
The general construction of such a weight system is very similar to that for Lie algebras, and we are not presenting it here.
Similarly to the case of a Lie algebra, computation of its values for a more or less complicated Lie superalgebra is rather laborious. 

In~\cite{FKV97}, this construction is elaborated for the simplest Lie superalgebra $\gl(1|1)$. The corresponding
weight system $w_{\gl(1|1)}$ takes values in the ring of polynomials in two variables, the subring in
the center $ZU\gl(1|1)$ of the universal enveloping algebra,
which is generated by the Casimir elements of grading~$1$ and~$2$.  

In~\cite{FKV97}, recurrence relations for computing
the values of this weight system, which are similar to the Chmutov--Varchenko relations for $w_{\sl(2)}$ are deduced. In~\cite{CL07}, it is proved that the values of this weight system depend on the intersection graph of the chord diagram only. The skew characteristic polynomial of graphs,
see Sec.~\ref{sssSCP}, is nothing but the extension of the
weight system $w_{\gl(1|1)}$ to a $4$-invariant of graphs.

\section{Embedded graphs and delta-matroids}\label{s4}

A chord diagram can be treated as an embedded graph 
with a single vertex on an orientable surface.
One of the key problems of the theory of weight systems
is how one can extend a given weight system to embedded
graphs with arbitrarily many vertices. Such
generalized weight systems are associated to 
finite type invariants of links in the same way as
ordinary weight systems are associated to finite type
invariants of knots.

There are various approaches to extending weight systems
and graph invariants to embedded graphs. The most
well developed approach consists in interpreting an
embedded graph as an ordinary graph to which information
about the embedding is added, see, e.g.~\cite{EMM13}.
Various extensions of the classical Tutte polynomial
has been constructed in this way~\cite{BR02,EMM15}.

In this section we describe another approach to constructing
extensions, which is under development in the last few years.
This approach is based on the analysis of the
Hopf algebra structure
on the space of delta-matroids. It leads to invariants
satisfying $4$-term relations and being generalized weight systems, therefore. Delta-matroids (also
known under the name of Lagrangian matroids) where introduced by Bouch\'et around 1990~\cite{B89,B91}. 
These combinatorial structures isolate equally well
essential features of both abstract and embedded graphs.
The Hopf algebra structure on various spaces of delta-matroids and $4$-term relations for them 
were introduced in~\cite{LZ17}.

\subsection{Embedded graphs and $4$-term relations}

\begin{figure}
    \centering
    \includegraphics[scale=0.4]{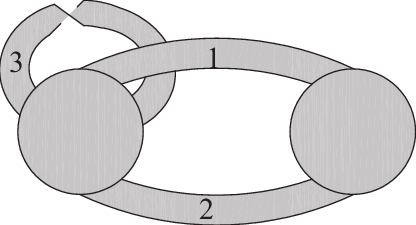}
    \caption{An embedded graph}
    \label{feg}
\end{figure}

Embedded graphs, also known under the name of ribbon graphs, are the subject of topological graph theory.
An \emph{orientable embedded graph} 
can be defined as an abstract graph (having, probably,
loops and multiple edges) endowed with a cyclic order
of semiedges incident to each vertex. 
Below, we consider only orientable connected embedded
graphs. Such a graph admits a presentation as an
orientable surface with boundary formed by disks 
associated to the verticies and ribbons associated to
the edges; each ribbon is attached to one or two
disks by its shorter sides, see Fig.~\ref{feg}.

To each chord diagram, one can associate a ribbon graph
with a single vertex by making the Wilson loop into
the vertex by attaching a disk to it and replacing
each chord with a ribbon in such a way that the surface remains orientable.

The $4$-term relations~(\ref{fourtermrelation}) 
can be interpreted as valid for an arbitrary embedded
graph, not necessarily a single-vertex one. 
In this more general case the ends of the two edges 
taking part in the relation may be attached to 
the boundaries of different vertices (two or even three).
We call these relations \emph{generalized $4$-term relations}, and we call the embedded graph invariants
satisfying them
\emph{generalized weight systems}.
Vassiliev's first move for embedded graphs exchanges
the two neighboring ends of two ribbons, and the second Vassiliev move consists in sliding an end of
a ribbon (a handle) along the second ribbon.

Below, we will be interested first of all in the question
how to extend a weight system to a generalized weight system.
The weight systems we are going to construct will be
based on the Hopf algebra structure on the space of delta-matroids.

\subsection{Delta-matroids of graphs and embedded graphs,\\ 
and $4$-term relations for them}

A \emph{set system} is a pair $(E;S)$
consisting of a finite set~$E$ and a subset $S\subset 2^E$ of the set of its subsets. A set system
 $(E;S)$ is said to be
\emph{proper} if the set~$S$ is nonempty;
below, we consider only proper set systems.

\begin{definition}
A proper set system $(E;S)$ is called a
\emph{delta-matroid} if the following
\emph{symmetric exchange axiom} is valid for it:

for any pair of subsets
$X,Y\in S$ and any element
$a\in X\Delta Y$ there is an element
$b\in Y\Delta X$ such that
$(X\Delta \{a,b\})\in S$.
\end{definition}

Here~$\Delta$ denotes the symmetric set difference operation, $X\Delta Y=(X\setminus Y)\cup(Y\setminus X)$. The element~$b$ in the symmetric exchange axiom 
can either coincide with~$a$ or be different from it. 
The elements of~$S$ are called \emph{admissible
sets} for the delta-matroid $(E;S)$.

To each simple graph~$G$ and any embedded graph~$\Gamma$,
one can associate a delta-matroid. For a graph,
the construction procceds as follows.
Recall that a graph~$G$ is said to be nondegenerate
if the determinant of its adjacency matrix over the two-element field~$\BF_2$ is~$1$. The \emph{delta-matroid
$\delta(G)=(V(G);S(G))$ of a graph~$G$} 
is the set system consisting of the set $V(G)$ 
of vertices of~$G$ and a set $S(G)\subset 2^{V(G)}$
of its subsets, where a subset $U\subset V(G)$ 
of the set of vertices of the graph is admissible
(i.e., $U\in S(G)$) iff the induced subgraph $G|_U$ is nondegenerate.

In turn, the  \emph{delta-matroid
$\delta(\Gamma)=(E(\Gamma);S(\Gamma))$ of an
embedded graph~$\Gamma$} is the set system consisting
of the set $E(\Gamma)$ of the edges of~$\Gamma$
and a set $S(\Gamma)\subset 2^{E(\Gamma)}$ 
of its subsets, where a subset $U\subset E(\Gamma)$ 
of the set of edges is admissible (i.e.,
$U\in S(\Gamma)$) iff the spanning embedded subgraph
$\Gamma|_U$ has a connected boundary.

If the genus of an embedded graph~$\Gamma$ 
is~$0$ (that is, this graph is embedded into the sphere),
a set of its edges is admissible iff the corresponding 
spanning subgraph $\Gamma|_U$ is a tree. For an arbitrary
genus, the admisssible subsets in the delta-matroid
$\delta(\Gamma)$ also are called
\emph{quasitrees}.

A.~Bouch\'et proved that the set systems $\delta(G)$  and $\delta(\Gamma)$
indeed are delta-matroids. In addition, if $\Gamma$ 
is an embedded graph with a single vertex, that is, a chord diagram,
the corresponding definitions are compatible:

\begin{theorem}
If~$C$ is a chord diagram, that is, an embedded graph
with a single vertex, then the delta-matroid~$\delta(C)$
is naturally isomorphic to the delta-matroid $\delta(g(C))$
of the intersection graph of~$C$.
\end{theorem}

In other words, the boundary of a chord diagram~$C$ is
connected iff its intersection graph is nondegenerate.
The last assertion has been rediscovered many times,
see, e.g., \cite{S01}.

The delta-matroid of a graph, as well as the delta-matroid of an orientable embedded graph is even.
A delta-matroid $(E;S)$ is said to be
\emph{even} if the number of elements in all its
admissible sets has the same parity,
$|X|\equiv|Y|~{\rm mod}~2$ for all $X,Y\in S$.
Since we are interested first of all in graphs
embedded into orientable surfaces, below, we talk
only about even delta-matroids, although most
of the constructions can be  as well extended to delta-matroids
that are not even, see Sec.~\ref{ssno}. 

Introduce the twist (or partial dual) of a delta-matroid
in the following way. Let $U\subset E$ be a subset
of the base set~$E$ of a delta-matroid
$D=(E;S)$. The \emph{twist} $D*U$ of a delta-matroid~$D$ in the 
\emph{subset}~$U$ is the set system
$D*U=(E;S*U)$, where $S*U=\{X\Delta U|X\in S\}$.
For delta-matroids of embedded graphs, the twist
operation corresponds to partial duality
introduced in~\cite{C09}. The twist of a delta-matroid
in an arbitrary subset of its base set is a delta-matroid~\cite{B91}.

\begin{definition}
The delta-matroid $\delta(G)$ of a graph~$G$ is called \emph{graphic}. The twist of a graphic delta-matroid
in a subset of its base set is called a 
\emph{binary delta-matroid}.
\end{definition}

In particular, the delta-matroid $\delta(\Gamma)$
of any embedded graph~$\Gamma$ on an orientable surface is even binary.

For delta-matroids, a $4$-term relation is defined.
Its definition requires defining two operations on
delta-matroids, namely, the first and the second Vassiliev moves.
Each of these operations is specified by a chosen
ordered pair of elements of the delta-matroid.

The result $\widetilde{D_{ab}}$ of
\emph{sliding the handle~$a$ along the handle~$b$}
in a delta-matroid~$D=(E;S)$, $a,b\in E$, 
is defined as $\widetilde{D_{ab}}=(E;\widetilde{S_{ab}})$,
where
$$ 
\widetilde{S_{ab}} =S\Delta\{X \sqcup a|X \sqcup b \in S 
\text{ and } X \subseteq E \setminus \{a, b\}\}.
$$
The result ${D'_{ab}}$ of
\emph{exchange of the ends of the handles~$a$ and~$b$}
in a delta-matroid~$D=(E;S)$, $a,b\in E$, 
is defined by ${D'_{ab}}=(E;{S'_{ab}})$,
where
$$ 
{S'_{ab}} =S\Delta\{X \sqcup \{a,b\}|X \in S 
\text{ and } X \subseteq E \setminus \{a, b\}\}.
$$
Both the first and the second Vassiliev move
take a binary delta-matroid to a delta-matroid of the 
same type.

A function~$f$ on even binary delta-matroids
\emph{satisfies $4$-term relations}
if for any such delta-matroid $D=(E;S)$
and any pair $a,b\in E$ of elements of the base set
we have
$$
f(D)-f(D'_{ab})=f(\widetilde{D_{ab}})-f(\widetilde{D'_{ab}}).
$$
If a function on even binary delta-matroids
satisfies $4$-term relations, then its restriction
to delta-matroids of embedded graphs on orientable
surfaces satisfies generalized $4$-term relations for them.

\subsection{Hopf algebras of delta-matroids}

Two set systems $(E_1;S_1)$ and $(E_2;S_2)$
are said to be \emph{isomorphic} if there is a one-to-one mapping of the set~$E_1$ to the set~$E_2$
that takes the set of subsets 
$S_1$ one-to-one to the set of subsets $S_2$. 

Denote by~$\cB^e_n$ the vector space spanned by
isomorphism classes of even binary delta-matroids
on $n$-element sets and define
$$
\cB^e=\cB^e_0\oplus\cB^e_1\oplus\cB^e_2\oplus\dots.
$$
The vector space~$\cB^e$ can be endowed with a
natural structure of graded commutative cocommuta\-tive
Hopf algebra. Multiplication in this algebra is
given by the disjoint union, while comultiplica\-tion
is given by splitting of the base set of a delta-matroid
into two disjoint subsets in all possible ways.
A delta-matroid is said to be \emph{connected}
if it cannot be represented as a product of two
delta-matroids with smaller base sets.

The quotient of the vector space $\cB^e$ modulo
the subspace of $4$-term relations inherits a
Hopf algebra structure. We are not going to use 
this Hopf algebra, and do not introduce a notation for it.

\begin{remark}
In contrast to graphs and delta-matroids, embedded
graphs, as far as we know, do not generate a natural
Hopf algebra.
\end{remark}

\subsection{Extending graph invariants to embedded
graphs and delta-matroids}

In this section, we explain how to extend to embedded
graphs the $4$-invariants of graphs described in Sec.~\ref{s2}.
The extension obtained happen to be defined on
binary delta-matroids, and their definition uses
the Hopf algebra structure on the vector space of binary
delta-matroids.

\subsubsection{Interlace polynomial}
The interlace polynomial of graphs defined above
admits a natural extension to binary delta-matroids.

\begin{definition}
The \emph{distance} from a set system $D=(E;S)$
to a subset $U\subset E$ is the number
$$
d_D(U)={\rm min}_{W\in S}~|U\Delta W|.
$$

The \emph{interlace polynomial} $L_D(x)$ 
of a delta-matroid~$D=(E;S)$ is the polynomial
$$
L_D(x)=\sum_{U\subset E}x^{d_D(U)}.
$$
\end{definition}

\begin{theorem}
If $D=D(G)$ is the delta-matroid of a graph~$G$,
then the interlace polynomial $L_D$ of~$D$ is 
related to the interlace polynomial of~$G$ in
the following way:
$$
L_{D(G)}(x)=L_G(x-1).
$$
\end{theorem}

\begin{theorem}[\cite{K20}]
The interlace polynomial of even binary delta-matroids
satisfies $4$-term relations for them.
\end{theorem}

\begin{theorem}
The space of values of the interlace polynomial on
the subspace of primitive elements in
$\cB^e_n$ has dimension~$n$.
\end{theorem}

This implies that the dimension of the space of
primitive elements in~$\cB^e_n$ is at least~$n$.

\subsubsection{Stanley's symmetrized chromatic polynomial}
The notion of umbral invariant of graphs admits a
natural extension to invariants of delta-matroids.
To define such an invariant, it suffices to specify,
for each connected delta-matroid~$D=(E;S)$ the coefficient
of~$q_n$, where $n=|E|$, in the value of the umbral
invariant of this delta-matroid.

In~\cite{NZ21}, a way to extending umbral graph invariants to umbral invariants of delta-matroids
has been suggested. This way is based on the combinatorial Hopf algebra structure introduced in~\cite{ABS06}.
Below, we restrict ourselves with commutative cocommutative Hopf algebras.

A \emph{combinatorial Hopf algebra} is a pair $(\cH,\xi)$ consisting of a connected graded
Hopf algebra~$\cH$ and a character (i.e.,
a multiplicative mapping) $\xi:\cH\to\BC$.

The polynomial Hopf algebra $\BC[q_1,q_2,\dots]$
becomes a combinatorial Hopf algebra if we endow it
with the \emph{canonical character} $\zeta$
which takes each variable~$q_k$ to~$1$.
This combinatorial Hopf algebra is a terminal object
in the category of combinatorial Hopf algebras, that
is, to each combinatorial Hopf algebra $(\cH,\xi)$
a unique umbral invariant (a graded Hopf algebra homomor\-phism) 
$\cH\to\BC[q_1,q_2,\dots]$ that takes~$\xi$
to~$\zeta$ is associated.

It is proved in~\cite{ABS06} that Stanley's symmetrized
chromatic polynomial is the graded Hopf algebra homomorphism $\cG\to\BC[q_1,q_2,\dots]$, which
corresponds to the character~$\xi:\cG\to\BC$,
which takes the graph with a single vertex to~$1$, 
and any connected graph with 
$n\ge2$ vertices to~$0$. 

This structure allowed one to define~\cite{NZ21} \emph{Stanley's symmetrized chromatic polynomial
of delta-matroids} as the graded homomorphism
$W:\cB^e\to\BC[x;q_1,q_2,\dots]$
from the Hopf algebra of even binary delta-matroids
to the Hopf algebra of polynomials in the variables~$q$
with coefficients in the ring of polynomials in a
variable~$x$, which is associated to the character
taking the delta-matroid $(\{1\};\{\{\emptyset\}\})$
to~$1$, the delta-matroid $(\{1\};\{\{1\}\})$
to~$x$, and any connected delta-matroid of grading
$n\ge2$ to~$0$. These two delta-matroids on 
$1$-element set span
the homogeneous subspace $\cB^e_1$.

\begin{theorem}
Stanley's symmetrized chromatic polynomial of
delta-matroids satisfies $4$-term relations for
delta-matroids.
\end{theorem}

\subsubsection{Transition polynomial}

Transition polynomial also admits an extension to
delta-matroids, which satisfies $4$-term relations for
them. In order to define transition polynomial, we will need few further definitions. Let $D=(E;S)$ be a 
delta-matroid, and let $U\subset E$ be a subset of
its base set.

The polynomial
$$
Q_D(s,t,x)=\sum_{E=\Phi\sqcup X\sqcup\Psi}
s^{|\Phi|}t^{|X|}x^{d_{D+\Phi*X\bar*\Psi}(\emptyset)}
$$
is called the \emph{transition polynomial}
of the delta-matroid
$D=(E;S)$. Here the operation~$+$, for an element $u\in E$, is defined by
$$
D+u=(E,S\Delta\{F\cup\{u\}|F\in S, u\notin F\}).
$$
The operation $\bar*$ is given by
$$
D\bar*u=D+u*u+u=D*u+u*u.
$$
These two operations are extended to an arbitrary subset $A\subset E$.

This definition is a specialization of the definition of weighted transition polynomial in~\cite{BH14}.
\begin{theorem}[~\cite{DZ22}]
The transition polynomial of even binary delta-matroids
satisfies $4$-term relations for them.
\end{theorem}

\subsubsection{Skew characteristic polynomial}

Graph nondegeneracy naturally extends to embedded graphs
and delta-matroids. A delta-matroid $(E;S)$ is said to be
\emph{nondegenerate} if $E\in S$, i.e., if the whole
set~$E$ is admissible. For embedded graphs, this means
that the embedded graph is nondegenerate iff its
boundary is connected. For chord diagram, this definition is consistent with the one we already know.
\emph{Nondegeneracy
$\nu(D)$ of a delta-matroid} is~$1$ if~$D$
is nondegenerate, and it is~$0$ otherwise.

We can define the \emph{skew characteristic 
polynomial of a delta-matroid} $D=(E;S)$ by
$$
Q_D(x)=\sum_{U\subset E} x^{|U|}\nu(D|_U).
$$

\begin{theorem}
The skew characteristic polynomial of delta-matroids
satisfies $4$-term relations.
\end{theorem}

Similarly to the case of graphs, this statement 
follows from the fact that nondegeneracy is
preserved under the second Vassiliev move.

\subsection{Extending to nonorientable surfaces:\\
framed chord diagrams, framed graphs,\\
and not necessarily even delta-matroids}\label{ssno}

$4$-term relations for chord diagrams, graphs, and embedded graphs have a nonorientable version, which
we describe here briefly. It should not be separated
from the orientable one: all the constructions, and first
of all the Hopf algebra structures, require simultaneous
participation of both orientable and nonorientable objects.
In each case, orientable objects generate Hopf subalgebras
in the Hopf algebras of objects without restrictions to their orientability.

A \emph{framed chord diagram} is a pair $(C,\varphi)$
consisting of a chord diagram~$C$ and a framing
$\varphi:V(C)\to\BF_2$ taking each chord of~$C$ to
an element of the field~$\BF_2$.
An isomorphism of framed chord diagrams is an
isomorphism of the underlying chord diagrams taking
the framing of the first of them to that of the second one.

A \emph{framed graph} is a pair $(G,\varphi)$ consisting of a simple graph~$G$ and a framing
$\varphi:V(G)\to\BF_2$, which associates to any
vertex of~$G$ an element of the field~$\BF_2$.
An isomorphism of framed graphs is an
isomorphism of the underlying graphs taking
the framing of the first of them to that of the second one.
By numbering the vertices of~$G$ with the numbers
$1,\dots,n=|V(G)|$, one associates to the framed graph its adjacency matrix
$A(G,\varphi)$. Nondiagonal elements of this matrix 
coincide with those of the adjacency matrix $A(G)$ of the simple graph~$G$, and its diagonal elements coincide
with the framings of the corresponding vertices.
Isomorphism classes of framed graphs coincide with
symmetric matrices over the field $\BF_2$
considered up to simultaneous renumbering of the
columns and the rows.

A framed graph is associated to any framed chord 
diagram, which is its intersection graph.
A framed chord diagram can be considered as a ribbon
graph with a single vertex: an ordinary ribbon is
associated to a chord with framing~$0$, while a chord
with framing~$1$ is replaced with a half-twisted ribbon. If there is a chord with framing~$1$, then the
surface corresponding to the framed chord diagram is 
nonorientable. 
A framed graph invariant~$f$ is called a $4$-\emph{invariant} if for any pair of vertices
$a,b\in V(G)$ it satisfies the $4$-term relation
$$
f((G,\varphi))-f((G'_{ab},\varphi))=
f((\widetilde{G}_{ab},\widetilde{\varphi}_{ab}))-
f((\widetilde{G}'_{ab},\widetilde{\varphi}_{ab})).
$$
There is a natural one-to-one correspondence between
the sets of vertices of any two graphs
participating in the relation, therefore,
both framings $\varphi$ and $\tilde\varphi_{ab}$ have a common
domain. Here $\widetilde{\varphi}_{ab}$ denotes the
framing of the graph $\widetilde{G}_{ab}$, which
coincides with $\varphi$ on all the vertices but~$a$,
and $\widetilde{\varphi}_{ab}(a)=\varphi(a)+\varphi(b)$.
In other words, the adjacency matrix
$A(\widetilde{G}_{ab})$ of the framed graph $\widetilde{G}_{ab}$ can be obtained from $A(G)$ by means of the same transformation~(\ref{eVsm})
as in the unframed case.

Nondegeneracy and skew characteristic polynomial
for framed graphs are defined in the same way as for
ordinary graphs. The nondegenracy of a graph is
preserved under the second Vassiliev move.
This implies that the skew characteristic polynomial
is a $4$-invariant of framed graphs. In contrast to
skew characteristic polynomial of simple graphs
the skew characteristic polynomial of a framed graph
can be neither even, nor odd.

Any $4$-invariant of framed graphs defines a framed
weight system, which takes a framed chord diagram
to the value of the invariant on its framed intersection graph.
In particular, nondegenercy and skew characteristic 
polynomial both define framed weight systems.

The delta-matroid of a framed graph and the delta-matroid
of an arbitary embedded graph, may be nonorientable,
are defined in the same way as in the orientable case.
If a framed graph has a vertex with framing~$1$,
then the corresponding delta-matroid will not be even:
its admissible sets will contain even as well as
odd number of elements. The delta-matroid of an 
embedded graph on a nonorientable surface also is
not even. Isomorphism classes of binary delta-matroids
generate a Hopf algebra containing the Hopf subalgebra
$\cB^e$.




\begin{thebibliography}{99}


\bibitem{Ab826} N.~H.~Abel, {\it Beweis eines Ausdruckes, von welchem die Binomialformel
   ein einzelner Fall ist},
	Journal f\"ur die reine und angewandte Mathematik {\bf 1} 159--160 (1826).

\bibitem{ABS06} M.~Aguiar, N.~Bergeron, F.~Sottile, {\it Combinatorial Hopf algebras
		and generalized Dehn–Sommerville relations}, Compositio Math. {\bf 142} 1--30 (2006).

\bibitem{AM13}

M.~Aguiar, S.Mahajan,
{\it Hopf monoids in the category of species},
in: HOPF ALGEBRAS AND TENSOR CATEGORIES,
Contemporary Mathematics, vol.~585, 17-- (2013)



\bibitem{ABS04}
R. Arratia, B. Bollob\'as, and G.B. Sorkin,
{\it A two-variable interlace polynomial},
Combinatorica,
24(4): 567--584, 2004


\bibitem{BN95}
   {D.~Bar-Natan},
  {\it On the Vassiliev knot invariants},
  Topology,
   {1995},
  {\bf 34},
   {423--472}


\bibitem{BNG96}
   {D.~Bar-Natan and S.~Garoufalidis},
  {\it On the Melvin--Morton–-Rozansky conjecture},
  Inventiones mathematicae,
   {1996},
  {\bf 125},
   {103-133}


\bibitem{BNV15}
   {D. Bar-Natan and H. T. Vo},
  {\it Proof of a conjecture of Kulakova et al. related to the $sl(2)$ weight system},
  Eur. J. Comb.,
   {2015},
  {\bf 45},
   {65--70}
   

\bibitem{B17}
A.~Bigeni,
{\it A generalization of the Kreweras triangle through the universal  $\mathfrak{sl}_2$ weight system},
Journal of Combinatorial Theory, Series A,
2017, vol. 161, 309--326

\bibitem{BM16} S.~Billey, P.~McNamara, 
{\it The contributions of Stanley to the fabric
		of symmetric and quasisymmetric functions}, in {\it The mathematical legacy of Richard P.
		Stanley}, Editors P.~Herch, T.~Lam, P.~Pylyavskyy, V.~Reiner  83--104 (2016).
		
		
\bibitem{BR02}
 B.~Bollob\'as, O.~Riordan,
{\it A polynomial of graphs on surfaces},
Mathematische Annalen,
vol. 323, 81--96 (2002)
		
\bibitem{B94}
   {A.~Bouchet},
  {\it Circle Graph Obstructions},
  J. Comb. Theory, Ser. B,
   {1994},
  {\bf 60},
   {107--144}

\bibitem{B91}
   {A.~Bouchet and A.~Duchamp},
  {\it Representability of delta-matroids over $GF(2)$},
  Linear Algebra and its Applications,
   {1991},
  {\bf 146},
   {67--78}

\bibitem{B89}
  {A.~Bouchet},
  {\it Maps and Delta-matroids},
   Discret. Math.,
   {1989},
  {\bf 78},
   {59--71}

\bibitem{BH14}
R. Brijder and H. Hoogeboom, 
{\it Interlace polynomials for multimatroids
and delta-matroids}, 
European Journal of Combinatorics, 
vol. 40, pp.
142--167, 2014.

\bibitem{BE11}
Viktor M Buchstaber and Nikolai Yu Erokhovets,
{\it Polytopes, Fibonacci numbers, Hopf algebras, and quasi-symmetric functions}, Russian Math. Surveys, 66:2 (2011), 271--367

\bibitem{BM19}
B. S. Bychkov, A. V. Mikhailov
{\it Polynomial graph invariants and linear hierarchies}, Russian Mathematical Surveys, 2019, Volume 74, Issue 2, Pages 366--368

\bibitem{CDL94}
   {S.~Chmutov, S.~Duzhin, and S.~Lando},
   {\it Vassiliev knot invariants. III: Forest algebra and weighted graphs},
Singularities and bifurcations, Adv. Soviet Math. 21 135--145,
Amer. Math. Soc., Providence, RI, 1994


\bibitem{CDBook12}
  {S.~Chmutov, S.~Duzhin, and J.~Mostovoy},
  {\it Introduction to Vassiliev Knot Invariants},
  Cambridge University Press,
    {2012}


\bibitem{CL07}
   {S.~Chmutov and S.~Lando},
  {\it Mutant knots and intersection graphs},
  Algebraic \& Geometric Topology,
   {2007},
  {\bf 7},
   {1579--1598}

\bibitem{C09}
S. Chmutov,
{\it Generalized duality for graphs on surfaces and the signed
Bollob\'as--Riordan polynomial},
J. of Combin. Theory Ser. B
{\bf 99} (2009)
617--638

\bibitem{CKL20}
S.~Chmutov, M.~Kazarian, and S.~Lando,
{\it Polynomial graph invariants and the KP
hierarchy},
Selecta Math. New Series
(2020)
26:34

\bibitem{CV97}
 S.~Chmutov, A.~Varchenko,
\newblock {\it Remarks on the Vassiliev knot invariants coming from $\mathfrak{sl}_2$}
\newblock Topology
\newblock 1997, 36, 1, 153--178



\bibitem{CMNR18}
   {C.~Chun, I.Moffatt, S. D. Noble, and R.~Rueckriemen},
  {\it On the interplay between embedded graphs and delta-matroids},
  Proceedings of the London Mathematical Society,
   {2018}


\bibitem{CMNR19}
   {C.~Chun, I.Moffatt, S. D. Noble, and R.~Rueckriemen},
  {\it Matroids, delta-matroids and embedded graphs},
  J. Comb. Theory, Ser. A,
   {2019},
  {\bf 167},
   {7--59}

\bibitem{DL22}
R. Dogra, S. Lando,
{\it Skew characteristic polynomial of graphs and embedded graphs}, arXiv: 2201.07084 (2022)


\bibitem{DZ22}
   {A.~ Dunaykin and V. Zhukov},
  {\it Transition polynomial as a weight system for binary delta-matroids},
  Moscow Math. J., vol.~22, no.~1, 1--13 (2022)

\bibitem{DBKPSS20}
 Petr Dunin-Barkowski, Maxim Kazarian, Aleksandr Popolitov, Sergey Shadrin, Alexey Sleptsov,
{\it Topological Recursion for the extended Ooguri-Vafa partition function of colored HOMFLY-PT polynomials of torus knots}
arXiv:2010.11021



\bibitem{EMM13}
J.~A.~Ellis-Monaghan, I.~Moffatt,
{\it Graphs on surfaces: dualities, polynomials, and knots}, 2013, Springer

\bibitem{EMM15}
J.~A.~Ellis-Monaghan, I.~Moffatt,
{\it The Las Vergnas polynomial for embedded graphs},
European Journal of Combinatorics
Volume 50, November 2015, Pages 97--114


\bibitem{FKV97}
Figueroa-O'Farrill, T.~Kimura, A.~Vaintrob,
\newblock {\it The Universal Vassiliev Invariant for the Lie Superalgebra $\gl(1|1)$}
\newblock Comm. Math. Phys.,
\newblock 1997, 185, 93--127

\bibitem{F20}
P.~Filippova,
\newblock {\it Values of the $\sl(2)$ Weight System on Complete Bipartite Graphs}
\newblock Funct. Anal. Appl.,
\newblock 54:3 (2020), 208--223

\bibitem{F22}
P.~Filippova,
\newblock {\it Values of the $\sl(2)$-weight system on a family of graphs that are not intersection graphs
of chord diagrams}
\newblock Sbornik Math.
\newblock 213:2 (2022),  115--148

\bibitem{HJ19}
S.~Heil and C.~Ji, 
{\it On an algorithm for comparing the chromatic symmetric functions
of trees},
Australas. J. Combin., 75: 210--222, 2019

\bibitem{H12}
J.~Huh,
{\it Milnor numbers of projective hypersurfaces and the chromatic polynomial of graphs},
J. Amer. Math. Soc. 25 (2012), no. 3, 907--927

\bibitem{J90} F.~Jaeger,
{\it On transition polynomial for $4$-regular graphs},
in: Cycles and rays (Montreal, PQ, 1987)m NATO Adv. Stud.
Inst. Ser. C: Math. Phys. Sci., 
Vol.~301, Kluwer Acad. Publ., Dordrecht 1990, 123--150 

\bibitem{JR79}
S.A.Joni and G.C.Rota,
{\it Coalgebras and Bialgebras in Combinatorics},
Studies in Applied Mathematics 61 (1979), 9--139


\bibitem{KP70}
B.~B.~Kadomtsev, V.~I.~Petviashvili,
\emph{On the stability of solitary waves in weakly dispersive media}, Sov. Phys. Dokl. 15: 539--541 (1970).

\bibitem{KL15} M.~Kazarian, S.~Lando, {\it Combinatorial solutions to integrable hierarchies}, Russ. Math. Surveys, {\bf 70}(3) 453--482 (2015).


\bibitem{K20}
N.~Kodaneva,
{\it The interlace polynomial of binary
delta-matroids and link invariants},
arXiv:2002.12440v1 (2020)

\bibitem{K93}
M.~Kontsevich.
\newblock {\it Vassiliev knot invariants},
\newblock in: {\em Advances in Soviet Math.}, 16(2):137--150, 1993.


\bibitem{K21}
E.~Krasilnikov,
{\it An Extension of the $sl_2$ Weight System to Graphs with $n\le 8$ Vertices},
Arnold Mathematical Journal, Vol.~7, No. 4, 609--618 (2021)



\bibitem{KLMR14}
   {E.~Kulakova, S.~Lando, T.~Mukhutdinova, and G. Rybnikov},
  {\it On a weight system conjecturally related to $sl_2$},
Eur. J. Comb.,
   {2014},
  {\bf 41},
   {266--277}


\bibitem{L06}
   {S.~Lando},
 {\it $J$-invariants of plane curves and framed chord diagrams},
    Functional Analysis and Its Applications,
   {2006},
  {\bf 40},
   {1--10}


\bibitem{L97}
  {S.~Lando},
  {\it On primitive elements in the bialgebra of chord diagrams},
   Translations of the American Mathematical Society-Series 2,
  {\bf 180},
   {167--174},
   {1997},
{Providence [etc.] American Mathematical Society, 1949-}


\bibitem{L00}
   {S.~Lando},
  {\it On a Hopf Algebra in Graph Theory},
  J. Comb. Theory, Ser.,
   {2000},
  {\bf 80},
   {104--121}


\bibitem{LZ03}
   {S.~Lando and A.~Zvonkin},
  {\it Graphs on Surfaces and Their Applications},
   {2003}


\bibitem{LZ17}
   {S.~Lando and V.~Zhukov},
  {\it Delta-matroids and Vassiliev invariants},
  Moscow Mathematical Journal,
  17(4) (2017)
741--755.

\bibitem{MM65}
J.~Milnor, J.~Moore,
{\it On the structure of Hopf algebras}, 
Ann. of Math. (2) 81 (1965), 
211--264

\bibitem{MMB17}
   {I.~Moffatt and E.~Mphako-Banda},
  {\it Handle slides for delta-matroids},
  Eur. J. Comb.,
   {2017},
  {\bf 59},
   {23--33}
   
\bibitem{M17}
A. Morse,
\emph{The Interlace Polynomial},
in: Graph POlynomials,
Eds. Y.Shi, M.Dehmer, X.Li, I.Gutman
Chapman and Hill, 2017

\bibitem{NZ21}
M.~Nenasheva, V.~Zhukov,
{\it An extension of Stanley's chromatic symmetric function to binary
delta-matroids}, Discete Math., 2021

\bibitem{NN}
N. Netrusova, 
{\it The interlace polynomial and knot invariants},
unpublished

\bibitem{NW99}
S. D. Noble; D. J. A. Welsh,
{\it A weighted graph polynomial from chromatic invariants of knots},
Annales de l'institut Fourier (1999),
Volume: 49, Issue: 3, page 1057--1087

\bibitem{OO96}
Okounkov, Andrei and Olshanski, Grigori
\newblock {\it Shifted Schur Functions},
\newblock June 1996.
\newblock St. Petersburg Math. J. volume 9.

\bibitem{O91}
G.~Olshanski,
\newblock {\it Representations of
infinite-dimensional classical groups,
limits of enveloping algebras and yangians},
\newblock in “Topics in Representation Theory”,
Advances in Soviet Math. 2,
\newblock Amer. Math. Soc., Providence RI, 1991, pp. 1--66.

\bibitem{RR78} S.~Roman, G.-C.~Rota, \emph{The Umbral Calculus}, Adv. Math. 27, 95--188~(1978)

\bibitem{RST97} G-C.~Rota, J.~Shen, B.~Taylor,
{\it All Polynomials of Binomial Type Are
    Represented by Abel Polynomials}, Annali della Scuola Normale Superiore di Pisa --
		Classe di Scienze, S\'erie 4 {\bf 25}(3–4): 731--738 (1997).


\bibitem{S83} M. Sato and Y. Sato, 
\emph{Soliton equations as dynamical systems on infinite
dimensional Grassmann manifolds}, in: Nonlinear partial differential equations
in applied science (Tokyo 1982), North-HollandMath. Stud., vol. 81, North-
Holland, Amsterdam 1983, pp. 259--271.

\bibitem{S94} W.~R.~Schmitt, 
{\it Incidence Hopf algebras},
		Journal of Pure and Applied Algebra {\bf 96} 299--330 (1994).

\bibitem{S95} W.~R.~Schmitt, 
{\it Hopf algebra methods in graph theory},
		Journal of Pure and Applied Algebra {\bf 101}(1) 77--90 (1995).


\bibitem{S01}
E. Soboleva,
{\it Vassiliev knot invariants coming from Lie algebras and $4$-
invariants},  J. Knot Theory and its Ramifications, 2001, {\bf 10}, 161--169.


\bibitem[St95]{St95} R.~Stanley, 
\emph{A symmetric function generalization of the chromatic
		polynomial of a graph}, Advances in Math. {\bf 111}(1) 166--194 (1995).

\bibitem{V90}
 V.~Vassiliev,
{\it Cohomology of knot spaces}
in: Advance in Soviet Math.
v.~1,
1990,
23--69


\bibitem{ZY21}
Zhuoke Yang,
\newblock {\it On values of $\sl_3$ weight system on chord diagrams whose intersection graph is complete bipartite},
\newblock arXiv:2102.00888

\bibitem{ZY22}
Zhuoke Yang,
\newblock {\it New approaches to $\gl(N)$ weight system}, 
arXiv:2202.12225 (2022)

\bibitem{Za22}
P.~Zakorko,
{\it Values of the $\sl(2)$ weight system on
chord diagrams with complete intersection graph}
\newblock to appear


\bibitem{Z70}
Zhelobenko D.P.
\newblock {\it Compact Lie groups and their representations}
\newblock Nauka, Moscow, 1970.
\newblock English translation: Translations of mathematical monographs, v. 40.
\newblock American Mathematical Society, Providence, Rhode Island, 1973.











\end{thebibliography}
\end{document}